\documentclass[conf]{new-aiaa}
\usepackage[utf8]{inputenc}
\usepackage{algorithm}
\usepackage{algpseudocode}
\usepackage{graphicx}
\usepackage{amsmath}
\usepackage[version=4]{mhchem}
\usepackage{siunitx}
\usepackage{longtable,tabularx}
\usepackage{multicol}
\usepackage{booktabs}
\usepackage{float}
\setcitestyle{authoryear,open={(},close={)}}
\usepackage{caption}
\usepackage{verbatim}
\usepackage{mathtools}

\usepackage{graphicx}
\usepackage{subcaption}
\usepackage{mwe}
\usepackage{xcolor}
\usepackage{lscape} 

\usepackage[mathscr]{euscript}

\defcitealias{EPA2018}{EPA, 2018}

\title{Electric Vehicle Traveling Salesman Problem with Drone with Partial recharge Policy}

\author{Tengkuo Zhu\footnote{Email: zhutengkuo@utexas.edu} }
 \affil{University of Texas at Austin, Texas.}
 \author{Stephen D. Boyles\footnote{Email: sboyles@mail.utexas.edu}}
 \affil{University of Texas at Austin, Texas}
 \author{Avinash Unnikrishnan\footnote{Email: droneinash@pdx.edu}}
 \affil{Portland State University, Portland, Oregon}

\begin{document}

\maketitle

\begin{abstract}
In \citep{zhu2022full}, it proposes an electric vehicle traveling salesman problem with drone while assuming that the electric vehicle (EV) is a battery-electric vehicle whose energy could be refreshed in a battery swap station in minutes. In this paper, we extend the work in \citep{zhu2022full} by relaxing the fixed-time-full-charge assumption, assuming that the EV is a plug-in hybrid electric vehicle that could be partially recharged in a charging station. This problem is named electric vehicle traveling salesman problem with drone with partial recharge policy (EVTSPD-P). A three-index MILP formulation is proposed to solve the EVTSPD-P with linear and non-linear charging functions where the concave time-state-of-charge (SoC) function is approximated using piecewise linear functions, a technique proposed in \cite{Montoya2017} and \cite{zuo2019new}. Furthermore, a specially designed adaptive large neighborhood search (ALNS) meta-heuristic, which incorporates constraint programming (CP), is presented to solve EVTSPD-P problem instances of practical size. The numerical analysis results indicate that the proposed ALNS method is more efficient than variable neighborhood search and has an average optimality gap of about 3\% when solving instances with ten nodes. Besides, using a piecewise linear function with a six-line-segments approximation has an average of 10.8\% less cost than a linear approximation.

\vspace{0.5cm}

\end{abstract}

\section{Introduction}
In the United States, the transportation sector generates 28.9\% of the national greenhouse gas emissions \citepalias{EPA2018}. Many local governments and corporate policies aim to promote transportation modes with lower pollution and greenhouse gas emissions. Electric vehicles (EVs) are an emerging alternative to internal combustion engines, and several companies have started to use EVs in their operations. For example, in 2018, FedEx announced a fleet expansion and added 1,000 electric delivery vehicles to operate commercial and residential pick-ups and delivery services in the United States \citep{FedEx}. Switching to electric fleets not only has long-term effects on mitigating the impact of climate change but may also have immediate financial benefits, as fuel cost accounts for 39\% to 60\% of operating costs in the trucking sector \citep{Sahin2009}. Compared to conventional internal combustion engines and petroleum-fuel-powered vehicles, EVs are much more energy-efficient and require less maintenance, which indicates potential savings to freight and logistics companies \citep{howey2011comparative,ma2012new}. 

Another new trend in recent years is the integration of UAVs into the operation of e-commerce and on-demand item delivery. The application of UAVs for "last-mile" parcel delivery promises to change the logistics industry landscape. Amazon, Google, and DHL all announced plans to use UAVs to deliver small packages, and Google has conducted thousands of test flights in Australia. The past few years have witnessed a dramatic increase in UAV applications \citep{Dronezon}. According to Teal Group's prediction, commercial use of UAVs will grow eightfold over the decade to reach US \$7.3 billion in 2027 \citep{Teal}. For the rest of the paper, the term "UAV" and "drone" are used interchangeably.

In \citep{zhu2022full}, it proposes EVTSPD, where an electric vehicle equipped with a drone is used to perform delivery tasks. Both the electric vehicle and the drone could serve customers independently, as proposed in \cite{Murray2015} and \cite{Agatz2018}. Furthermore, EVTSPD assumes that the EV shares the energy with the drone; that is, the onboard battery is the driving force for the EV and can provide energy for the drone. Besides, the EV could recharge at the charging station nodes en route when necessary. In \citep{zhu2022full}, it focuses on EVTSPD with a fixed-time-full-charge policy, which assumes that the EV could refresh its battery to full energy with fixed amount of time. In this paper, we extend the work in \citep{zhu2022full} by relaxing the fixed-time-full-charge assumption. Instead, this paper assumes that the EV can be partially charged with linear/non-linear time-state-of-charge functions while the charging time is also a decision variable that depends on the state-of-charge level and charging functions. This problem is named electric vehicle traveling salesman problem with drone with partial recharge policy, or EVTSPD-P. 

In EVTSPD-P, the electric vehicle used in EVTSPD instances is considered a plug-in hybrid electric vehicle, and the charging station nodes in the network are assumed to be connected to a plug-in electric grid. Thus, the electric vehicle in EVTSPD-P can be charged at CS nodes with any energy, while the required charging time depends on the initial state-of-charge level before charging and the amount of charged energy. This realistic assumption indicates that the proposed EVTSPD-P model can be applied in real-world scenarios. This paper aims to shed some light on how several key factors, such as charging function, charging policy or EV's driving range, etc., affect the overall performance of the delivery system. Specifically, the main contributions of this paper are: 
\begin{itemize}
    \item In EVTSPD-P, we consider two different charging functions. The first one models linear time-SoC function while the other one models non-linear function. These two problems are named EVTSPD-PL and EVTSPD-PN, respectively.
    \item An alternative MILP formulation for EVTSPD, defined in an augmented network, is proposed. In the augmented network, the charging station nodes are copied to enable potential multiple visits to each CS node. This formulation enables us to incorporate the partial recharge assumption with linear/non-linear time-SoC function and is an essential component when modeling EVTSPD-PL and EVTSPD-PN.
    \item As EVTSPD-PL and EVTSPD-PN are complicated in a strong NP-hard sense, commercial solvers can only handle instances of small size. Thus, a specially designed adaptive large neighborhood search method is proposed, which incorporates constraint programming modeling as a subproblem. The test results indicate that this approach is more efficient than the variable neighborhood search method, which is widely used to solve traveling salesman problem with drone (TSPD).
    \item Additional numerical analysis shows that with the proposed MILP model, the commercial solver can solve instances containing up to 10 nodes in the network. Furthermore, parameter analysis shows that using piecewise linear approximation can obtain a solution whose cost is about 10.8\% less than linear time-SoC function approximation.
\end{itemize}\par
The rest of the paper is organized as follows: Section 2 describes EVTSPD-PL and EVTSPD-PN and introduces their MILP formulation for both models. Different time-SoC function approximation techniques are also explained, along with the additional efforts needed to incorporate them in the proposed MILP formulation.
Section 3 describes the proposed ALNS method. Section 4 presents the numerical analysis result, which includes the performance evaluation of the MILP formulation and ALNS algorithm. It also includes the sensitivity analysis of several key parameters used in the model. Section 5 presents the conclusion and the potential future research directions for the problem.\par

\section{Literature Review}
The vehicle routing problem (VRP) and the traveling salesman problem (TSP) are among the most well-studied optimization problems in operations research. This problem was first proposed by \citep{Dantzig1959}. Since then, many variants have been considered, incorporating service time windows, capacities, maximum route lengths, distinguishing pick-ups and deliveries, fleet inhomogeneities, etc. Various exact and heuristic methods have been proposed to solve the problem  \citep{baldacci2012recent, laporte2009fifty, desaulniers2010vehicle, Osman1993, Gendreau1992}. \citep{braekers2016vehicle, montoya2015literature, pillac2013review, toth2014vehicle, eksioglu2009vehicle, golden2008vehicle, toth2002vehicle} provide a thorough literature review of VRP variants and solution algorithm families.

\citep{Erdogan2012} introduced the green vehicle routing problem (G-VRP), where the goal is to route a fleet of Alternative Fueled Vehicles (AFV) to serve a set of customers within a time limit while respecting the driving range of the vehicles. The AFVs are allowed to extend their driving range by visiting refueling stations more than once.  \citep{conrad2011recharging} developed the recharging vehicle routing problem where vehicles can recharge at particular customer locations. \citep{conrad2011recharging} also considered customer time windows and fleet capacity constraints. \citep{schneider2014electric} also modeled customer time windows and fleet capacity constraints in Electric VRP (E-VRP) problem while making the recharging times dependent on the remaining charge levels. \citep{Montoya2017} extended the previous E-VRP models to consider nonlinear charging functions by using piecewise linear approximations. The solution methods for E-VRP variants are diverse and range from exact methods such as branch and bound, branch and cut \citep{kocc2016green}, and branch and price \citep{schneider2014electric, hiermann2016electric}; to heuristic methods such as modified savings method of Clarke and Wright with density-based clustering \citep{Erdogan2012}, local improvement based on neighborhood swap \citep{schneider2014electric, masmoudi2018dial}, and metaheuristics such as simulated annealing and tabu search \citep{keskin2016partial, goeke2015routing, felipe2014heuristic}.  \citep{Pelletier2017} and \citep{erdelic2019survey} provide a comprehensive survey of the different variants of the electric vehicle routing problem and associated solution algorithms. \emph{Unlike the research mentioned above, this paper focuses on the traveling salesman variant}. \citep{doppstadt2016hybrid} formulated the traveling salesman problem for hybrid electric vehicles considering four modes of operation - combustion, electric, charging, and boost. An iterated tabu search with local search operators, which switch route structure and operating modes, was used to solve real-world instances. \citep{doppstadt2019hybrid} extend \citep{doppstadt2016hybrid}'s model by considering customer time windows and proposed a new variable-neighborhood-search-based solution method. \citep{liao2016electric} provided an efficient dynamic-programming-based polynomial-time algorithm for the electric vehicle shortest travel time path problem and approximation algorithms for the EV touring problems. The algorithms incorporated battery capacity constraints and battery swaps. \citep{roberti2016electric} provided a mixed-integer linear programming formulation for the electric vehicle traveling salesman problem with time windows for both full and partial recharge policies. A three-phase heuristic employing variable neighborhood descent to reach time window feasibility and minimize cost tour and a dynamic programming algorithm to achieve feasibility concerning battery capacities is developed. \emph{While there has been a significant amount of work on E-VRP and TSP, except for \citep{zhu2022full}, none of them have considered an integrated delivery system with drones.}

Meanwhile, an increasing number of studies investigate the efficiency of delivery systems that deploy UAVs. \citep{otto2018optimization} provides a detailed review of the various civil applications of drones in domains such as agriculture, monitoring, transport, security, etc.  \citep{Murray2015} introduced the flying sidekick traveling salesman problem (FSTSP), which assumes that a truck can launch its UAV at the depot or customer node and remains on its route while the UAV delivers one small parcel to another customer before meeting again at a rendezvous location (another customer node on the truck's route). \citep{Murray2015} proposed a two-stage route and reassign heuristic wherein the first stage, a truck TSP tour that visits all customers, is determined. In the second stage, selected customers are reassigned to the UAV based on cost savings.  \citep{Murray2015} also introduced the parallel drone scheduling Traveling Salesman Problem (PDSTSP), where multiple drones and a truck originating from a depot serve a set of customers. In the PDSTSP heuristic, customers are partitioned into those that can be served by UAVs and those assumed to be served by the truck. A parallel machine scheduling problem is solved to determine customer assignments to drones. A swap-based heuristic is used to exchange customers from UAV and truck partitions to improve the solution. \citep{mbiadou2018iterative} developed an iterative two-stage heuristic involving customer partitioning and routing optimization for the PDSTSP.  \citep{Agatz2018} developed an integer programming formulation for a variant of FSTSP called Traveling Salesman Problem with Drones (TSPD) and a "route-first, cluster-second" heuristic, which constructs a TSP with drone tour from a TSP tour. A subtle difference between FSTSP and TSPD is that in FSTSP, the drone departs from a truck at a node and joins the truck at a different node, whereas in TSPD, the truck can wait at a node, and the drone can rejoin the truck at the same node it departed from. Note that several authors have used TSPD while referring to FSTSP. \citep{Ha2018} studied a variant of \citep{Murray2015} with the objective of minimizing operating and waiting time costs rather than completion time. The authors propose two heuristics - a modification of \citep{Murray2015}'s heuristic to minimize costs and Greedy Randomized Adaptive Search Procedure (GRASP). \citep{EsYurek2018} developed an iterative two-stage algorithm to solve \citep{Murray2015}'s FSTSP, which was referred to as TSPD. In the first stage, the truck route is determined, whereas, in the second stage, the drone tours are determined. \citep{bouman2018dynamic} modify the Bellman‐Held‐Karp dynamic programming algorithm for the TSP to develop an exact solution approach for TSPD, whereas \citep{poikonen2019branch} use a branch and bound method. \citep{de2020variable, DeFreitas2018} developed a randomized variable neighborhood descent heuristic, which modified an initial TSP solution obtained from the Concorde solver to solve the FSTSP. \citep{boysen2018drone} focus on scheduling single and multiple drone deliveries launched from a truck with a fixed route. \citep{jeong2019truck} modify \citep{Murray2015}'s model to include the impact of payload on energy consumption and no-fly zones and propose a two-stage construction and search heuristic. \citep{dayarian2020same} focused on a new variant where drones are used to resupply a truck making deliveries. \citep{kim2018traveling} study a variant termed Traveling Salesman Problem with Drone Station (TSPDS), where a truck is used to resupply a drone station, which is different from a depot. Multiple drones then make deliveries to customers from drone stations. The truck will also make deliveries to customers after supplying the drone station. A two-phase solution algorithm is developed, which involves determining optimal TSP for customers who can be served by truck only and parallel machine scheduling problem to determine customer assignment to drones. \emph{While there has been a significant body of work on integrating drones into existing routing frameworks since 2015, none of them consider EVs and their associated range constraints.}

Other researchers have used continuous approximation techniques to analyze UAV routing problems. \citep{carlsson2018coordinated} used a continuous approximation approach to study improvements in efficiency by using a drone with a traveling salesman problem framework. Based on asymptotic as well as computational analysis, the efficiency improvements were found to be proportional to the square root of the ratio of the speeds of the truck and the UAV. \citep{ferrandez2016optimization} determined that using multiple drones per truck led to an increase in savings in energy and time and developed continuous approximation formulas to estimate the savings. \citep{figliozzi2017lifecycle} compared lifecycle $CO_2$ emissions of drones relative to other delivery mechanisms such as diesel vans and electric trucks. UAVs were found to have lower lifecycle $CO_2$ emissions per distance compared to typical diesel vans.

Several researchers have focused on VRP variants involving drones. \citep{dorling2016vehicle} formulated drone delivery problems as a multi-trip VRP. Key contributions were a linear approximation of energy consumption as a function of payload and a simulated-annealing-based solution algorithm. \citep{Wang2017} derive worst-case bounds on the maximum savings obtained by integrating drones into traditional truck deliveries. \citep{ham2018integrated} adopt a constraint programming approach where multiple drones and trucks depart from a single depot. Drones can deliver as well as pick up while considering customer time windows. \citep{ulmer2018same} model deliveries using a heterogeneous fleet of trucks and drones from a single depot and found that spatial partitioning of the delivery zones into those delivered by trucks and those delivered by drones are more effective. \citep{wang2019vehicle} studied a VRP with a drone variant where drones can be launched from a truck, serve multiple customers, and then return to a docking hub from which they can be picked up by the same or different trucks. A branch-and-price formulation is developed to solve the mixed-integer linear program. \citep{sacramento2019adaptive} formulated the VRP variant of FSTSP where multiple trucks and UAV combinations are used to serve customers and solved the model using an adaptive large neighborhood search metaheuristic. 

In this paper, we extend the work in \citep{zhu2022full} by assuming that the EV could be partially charged at the charging station nodes. It turns out that the partial-recharge assumption greatly increases the difficulty of the problem as additional decision variables are needed to model the linear/nonlinear time-state-of-charge function. This assumption also indicates that the model proposed in this paper has great potential to be used in a real-world application.

\section{Problem Description and Formulations of EVTSPD-P}
As proposed in \citep{zhu2022full}, EVTSPD, which involves a single EV and a single drone, aims to find a coordinated route such that each customer is visited by the EV or the drone exactly once. The operation process of EVTSPD is similar to TSPD, where both the EV and the drone could serve the customer while the EV serves as the drone hub that can launch and retrieve the drone. Other assumptions adopted in EVTSPD include:
\begin{itemize}
    \item The EV could recharge its battery at the charging station nodes in the network. This charging process is assumed to start immediately for EV upon arrival.
    \item The EV and the drone are supposed to travel in constant speed, and the altitude differences in the network are ignored.
    \item Whenever the EV travels across an arc in the network, its battery energy decreases gradually. Besides, the total consumed energy of traversing an arc is a known constant value.
    \item Initially, both the EV and the drone are located at the depot and (only) the EV is assumed to be fully charged.
\end{itemize}

Besides, one fundamental assumption adopted in EVTSPD is the shared energy assumption, which indicates that the battery on EV also recharges the drone if it is on board. During the operation, whenever the drone is launched to serve a customer independently, the energy consumed by the drone should be deducted from the remaining energy of the battery. More explanations of this assumption can be seen in \citep{zhu2022full}.

In \citep{zhu2022full}, it focuses on EVTSPD-FF, which is a variant of EVTSPD and assumes that the EV could refresh its battery energy to full capacity within a fixed amount of time whenever it visits the charging station nodes. A MILP formulation is proposed for EVTSPD-FF, where the key idea is to create a multigraph with no charging station nodes. However, in this paper, as we discard the fixed-time-full-charge assumption and instead assume that the EV could be partially charged at the charging station nodes, the multigraph approach is no longer applicable in EVTSPD-P. Thus, this section first presents an alternative MILP formulation for EVTSPD and then discusses incorporating the partial recharge assumptions into this formulation. Note that this formulation is still valid for EVTPSD with a full-charge policy.

\subsection{MILP formulation for EVTSPD}
This section presents a three-index MILP formulation for EVTSPD,  which is defined in an augmented network. To permit multiple (and possibly zero) visits to the charging station nodes while requiring exactly one visit to the customer nodes, the original network $G$ is augmented to create $G^{a} = (N, E^{a})$. In this network, denote $N = V \cup S^{a}$ as the augmented vertices set where $S^{a}$ is the augmented set of all charging stations, including the $m$ copies of set $S$. Each copy indicates a potential visit to each charging station node. The number of dummy vertices associated with each charging station, $m$, is set to the number of times the associated $s \in S$ charging station node can be visited. $m$ should be set as small as possible to reduce the augmented network size but large enough not to restrict multiple beneficial visits. In the past literature, \cite{Erdogan2012}, this value is typically set to 2 or 3. Additional notations used in this model are defined next.

The sets used in the three-index EVTSPD model are defined below:
{\renewcommand\arraystretch{1.0}
\noindent\begin{longtable*}
{@{}l @{\quad:\quad} p{15cm}@{}}
% {@{}>{\raggedright}p{2cm}lp{7.3cm}@{}}
$C$ & Set of all customers in the problem, $C = \{1, 2, ... ,c\}$ and $|C| = c$\\
$C^{'}$ & Subset of customers that are available to drone delivery service, $C^{'} \subset C$\\

$S$ & Set of all charging stations in the original network, $S = \{c+1,c+2, . . . ,c+s\}$ and  $|S| = s$\\
$S^{a}$ & Augmented set of all charging stations which includes a total of $m$ copies of set $S$. $S^{a} = \{c+1,.., c+s, c+s+1, ..., c+ms \}$ and $|S^{a}| = ms$. Note that in the augmented network each node can be visited at most once \\
$N$ & Set of all nodes in the augmented network, $N = S^{'} \cup C \cup \{0, 0^{'}\}$, where $0$ and $0^{'}$ both represent the depot and $|N| = c+ms+2$\\
$N_{d}$ & Set of nodes from which a vehicle may depart in the augmented network. $N_{d}=\{0, 1, . . . ,c\} \cup S^{'}$\\
$N_{a}$ & Set of nodes to which a vehicle may arrive in the augmented network. $N_{a}=\{1, 2,...c\} \cup S^{'} \cup \{ v_{c+ms+1}\}$\\
$N^{'}$ & Set of all customer nodes and charging station nodes in the augmented network. $N^{'}= C \cup S^{a}$\\
$D$ & Set of drone's feasible sortie $\textlangle i,j,k\textrangle$, where the drone is launched from node $i$, travels to node $j$ and returns to node $k$. $D = \{ \textlangle i,j,k \textrangle : i \in N_{d}, j \in C^{'}, k \in N_{a}, i \neq j, j \neq k, i\neq k,e^{D}_{ij} + e^{D}_{jk} \leq Q^{D} \}$, where $Q^{D}$ represents flight energy limit of the drone and $e^{D}_{ij}$ represents the drone's energy cost when flying from node $i$ to node $j$. 
\end{longtable*}}

Below is a summary of parameters that are used in the modeling process:
{\renewcommand\arraystretch{1.0}
\noindent\begin{longtable*}
{@{}l @{\quad:\quad} p{15cm}@{}}
$c^{T}_{ij}$ & Time cost of EV when travelling from node $i \in N$ to node $j \in N$\\
$c^{D}_{ij}$ &  Time cost for the drone when travelling from node $i$ to node $j$ \\
$d_{ijk}$ &  Travel time cost for the drone to launch from node $i$, serves node $j$ and return to node $k$\\
$s_{i}$ &  Service time for customer $i \in C$\\
$e^{T}_{ij}$ & Energy cost for the EV when travelling from node $i$ to node $j$\\
$e^{D}_{ij}$ & Energy cost for the drone when travelling from node $i$ to node $j$\\
$e^{D}_{ijk}$ &  Energy cost for the drone to accomplish sortie $\textlangle i,j,k \textrangle \in D$\\
$Q^{D}$  & Flight energy capacity of the drone \\
$Q^{T}$  & Driving energy capacity of the EV \\ 
$M$ & A positive large number\\
\end{longtable*}}

The decision variables of the formulation is listed below:
\noindent\begin{longtable*}
% {@{}l @{\quad:\quad} l@{}}
{@{}l @{\quad:\quad} p{15cm}@{}}
$x_{ij} \in \{0,1\}$  & Equals one if the EV travels from node $i$ to node $j$ and zero otherwise, where $i \neq j$ and $i \in N_{d}, j \in N_{a}$ \\
$y_{ijk} \in \{0,1\}$  & Equals one if the drone performs sortie $\textlangle i,j,k \textrangle \in D$ and zero otherwise \\
$p_{ij}\in \{0,1\}$ & Equals one if customer node $i$ is visited  before customer $j$ in the EV's path and zero otherwise\\
$u_{i}$ & Position of node $i$ in the EV's path\\
$b_{i}^{a} \geq 0$ & Remaining battery charge of the EV upon arrival at node $i$ \\
$b_{i}^{d1} \geq 0$ & Remaining battery charge of the EV upon departure from node $i$ before launching drone\\
$b_{i}^{d2} \geq 0$ & Remaining battery charge of the EV upon departure from node $i$ after launching drone\\
$t^{a}_{j} \geq 0$ & Time when the EV arrives at node $j$ \\
$t^{d}_{j} \geq 0$ & Time when the EV departs from node $j$ \\
$t_{j}^{'} \geq 0$ & Time when the drone arrives at node $j$ 
\end{longtable*}

Note that parameter and decision variables associated with driving/flight energy, such as $b_{i}^{a}$, $b_{i}^{d1}$ or $e^{T}_{ij}$, have same measurement units. This paper assumes that the EV and drone have constant traveling time. The driving/flight energy of the EV/drone associated with a node, or the driving/flight energy cost associated with each arc, can be measured in the same time units as $c^{T}_{ij}$ and $c^{D}_{ij}$.

\newpage
Now we are ready to present the three-index MILP formulation. 
% Note that this formulation already considers fixed charging/service time, incompatible customers, and launch/retrieve times and does NOT consider loops in its formulation.

\noindent \textbf{\underline{Three-index model}}:

\vspace{0.2cm} 
\noindent \textbf{Objective:}
\begin{align}
\min \quad
& t_{c+ms+1}
\end{align}

\vspace{0.5cm}
\noindent The objective is to minimize the time when both EV and drone return to the depot after serving all the customers. There are many constraints in the problem, which are introduced below and interspersed with descriptions. \par

\hfill \break
\noindent \textbf{Routing Constraints:}
\begin{align}
\sum_{\substack{i \in N_{d}\\ i \neq j}} x_{ij}+\sum_{\substack{i \in N_{d}\\ i \neq j}} \sum_{\substack{k \in N_{a}\\\textlangle i,j,k \textrangle \in D} } y_{ijk} &= 1 && \forall j \in C \label{eqn:secmodel_con1}\\
%------------------------------------%-----------------------------------------%
\sum_{j\in N_{a}} x_{0j}& = 1  \label{eqn:secmodel_con2} \\
%------------------------------------%-----------------------------------------%
\sum_{i\in N_{d}} x_{i,0^{'}}& = 1 \label{eqn:secmodel_con3} \\
%------------------------------------%-----------------------------------------%
u_{i} - u_{j} +1 & \leq (c+ms+2)(1-x_{ij}) && \forall i \in N^{'}, j \in N_{a}, j \neq i \label{eqn:secmodel_con4} \\
%------------------------------------%-----------------------------------------%
\sum_{\substack{i\in N_{d}\\ i \neq j}} x_{ij} & = \sum_{\substack{k\in N_{a}\\ k \neq j}} x_{jk} && \forall j \in N^{'} \label{eqn:secmodel_con5} \\
%------------------------------------%-----------------------------------------%
\sum_{j \in C^{'}} \sum_{\substack{k \in N_{a} \\ \textlangle i,j,k \textrangle \in D} } y_{ijk} & \leq 1 && \forall i \in N_{d} \label{eqn:secmodel_con6} \\
%------------------------------------%-----------------------------------------%
\sum_{i \in N_{d}} \sum_{\substack{j \in C^{'} \\ \textlangle i,j,k \textrangle \in D} } y_{ijk} & \leq 1 && \forall k \in  N_{a} \label{eqn:secmodel_con7}\\
%------------------------------------%-----------------------------------------%
2 y_{ijk} & \leq \sum_{\substack{h \in N_{d}\\ h \neq i}} x_{hi} + \sum_{\substack{l \in N_{d}\\ l \neq k}} x_{lk} && \forall i \in N_{d}, j \in C^{'}, k \in N_{a}, \textlangle i,j,k \textrangle \in D \label{eqn:secmodel_con8}\\
%------------------------------------%-----------------------------------------%
y_{0jk} & \leq \sum_{\substack{h \in N_{d}\\ h \neq k}} x_{hk}  && \forall j \in C^{'}, k \in N_{a}, \textlangle 0,j,k \textrangle \in D \label{eqn:secmodel_con9}\\
%------------------------------------%-----------------------------------------%
1-(c+ms+2) \left( 1- \sum_{\substack{j \in C^{'}\\\textlangle i,j,k \textrangle \in D} } y_{ijk}\right) & \leq u_{k} - u_{i} && \forall  i \in N_{d}, k\in N_{a}, k \neq i \label{eqn:secmodel_con10}
\end{align}

Constraints (\ref{eqn:secmodel_con1})--(\ref{eqn:secmodel_con10}) are associated with the routing of the two vehicles. In particular, constraint (\ref{eqn:secmodel_con1}) guarantees that each customer node is visited once by either the EV or the drone. Constraints (\ref{eqn:secmodel_con2}) and (\ref{eqn:secmodel_con3}) state that the EV must start from and return to the depot.  Constraint (\ref{eqn:secmodel_con4}) is a sub-tour elimination constraint for the EV. Constraint (\ref{eqn:secmodel_con5}) indicates that if the EV visits node $j$ then it must also depart from node $j$. Constraints (\ref{eqn:secmodel_con6}) and (\ref{eqn:secmodel_con7}) state that each node can launch or retrieve the drone at most once. Constraint (\ref{eqn:secmodel_con8}) ensures that if there exists a drone route $\textlangle i,j,k \textrangle$, then EV should travel between $i$ and $k$. Constraint (\ref{eqn:secmodel_con9}) states that if the drone is launched from the depot and returned to node $k$, then node $k$ should be visited by the EV. Constraint (\ref{eqn:secmodel_con10}) is a sub-tour elimination constraint for the drone.\par
\hfill \break
\noindent \textbf{Energy Constraints:}
\begin{align}
% Battery level constraint start
b^{a}_{j} & \leq b^{d2}_{i} - e^{T}_{ij} x_{ij} + M (1-x_{ij}) && \forall i \in N_{d}, j \in N_{a}, i \neq j \label{eqn:secmodel_con11}\\
%------------------------------------%-----------------------------------------%
b^{a}_{0} & = Q^{T} \label{eqn:secmodel_con12}\\
%------------------------------------%-----------------------------------------%
b^{a}_{i} & = b^{d1}_{i} && \forall i \in C_{0} \label{eqn:secmodel_con13} \\
%------------------------------------%-----------------------------------------%
b^{a}_{i} & \leq b^{d1}_{i} && \forall i \in S \label{eqn:secmodel_con14} \\
%------------------------------------%-----------------------------------------%
b^{d2}_{i} & = b^{d1}_{i} - \gamma * \sum_{\substack{j \in C^{'} \\ j \neq i}}\sum_{\substack{k \in N_{a} \\ \textlangle i,j,k \textrangle \in D}}y_{ijk}e^{D}_{ijk} && \forall i \in N_{d} \label{eqn:secmodel_con15}\\
%------------------------------------%-----------------------------------------%
0 \leq b^{a}_{i}, b^{d1}_{i}, b^{d2}_{i} & \leq Q^{T}  && \forall i \in N \label{eqn:secmodel_con16}
\end{align}

Constraints (\ref{eqn:secmodel_con11})--(\ref{eqn:secmodel_con16}) are associated with the battery energy level. In particular,  constraint (\ref{eqn:secmodel_con11}) states that if the EV travels from node $i$ to node $j$, then the energy level before arriving at node $j$ is $\tau_{ij}$ less than the energy level after leaving node $i$, regardless whether node $i,j$ are customer nodes or charging stations. Constraint (\ref{eqn:secmodel_con12}) ensures that when EV departs from the depot it is fully charged. Constraints (\ref{eqn:secmodel_con13}) indicate that if the EV visits a customer node $i$, the remaining energy upon departure at node $i$ is the same the energy upon arrival at node $i$. Similarly, constraints (\ref{eqn:secmodel_con14}) enforce that if the EV visits a charging station $i$, then the remaining energy upon departure at node $i$ is no less than the energy upon arrival at node $i$. Constraints (\ref{eqn:secmodel_con15}) models the energy change at each node $i \in N_{d}$ because of the launch of the drone. Constraints (\ref{eqn:secmodel_con16}) are the domain constraints for $b^{a}_{i}, b^{d1}_{i}, b^{d2}_{i}$.

% Constraint (\ref{eqn:secmodel_con13}) states that if the EV departs from a charging station node $i$ and there is a drone route that starts at node $i$, then when EV departs from node $i$ it is no longer fully charged and the drone route energy consumption should be deducted from full-charged battery. Constraint (\ref{eqn:secmodel_con15}) states the same situation as constraint (\ref{eqn:secmodel_con14}) except when node $i$ is a customer. Constraint  (\ref{eqn:secmodel_con16}) ensures that the remaining battery charge should be non-negative and no greater than $Q^{T}$.\par

\hfill \break
\noindent \textbf{Coordination Constraints:}
\begin{align}
t^{'}_{i} &\geq t^{d}_{i}- M \left(1-\sum_{\substack{j \in C^{'} \\ j \neq i}} \sum_{\substack{k \in N_{a}\\\textlangle i,j,k \textrangle \in D} } y_{ijk}\right) && \forall i \in N_{d} \label{eqn:secmodel_con17}\\
%------------------------------------%-----------------------------------------%
t^{'}_{i} &\leq t^{d}_{i}+ M \left(1-\sum_{\substack{j \in C^{'} \\ j \neq i}} \sum_{\substack{k \in N_{a}\\ \textlangle i,j,k \textrangle \in D} } y_{ijk}\right) && \forall i \in N_{d}\label{eqn:secmodel_con18}\\
%------------------------------------%-----------------------------------------%
t^{'}_{k} &\geq t^{a}_{k}- M \left(1-\sum_{\substack{i \in N_{d} \\ i \neq k}} \sum_{\substack{j \in C^{'} \\ \textlangle i,j,k \textrangle \in D} } y_{ijk}\right) && \forall k \in N_{a} \label{eqn:secmodel_con19}\\
%------------------------------------%-----------------------------------------%
t^{'}_{k} &\leq t^{a}_{k}+ M \left(1-\sum_{\substack{i \in N_{d} \\ i \neq k}} \sum_{\substack{j \in C^{'} \\ \textlangle i,j,k \textrangle \in D} } y_{ijk}\right) && \forall k \in N_{a} \label{eqn:secmodel_con20}\\
%------------------------------------%-----------------------------------------%
t^{a}_{k} &\geq t^{d}_{h}+c^{T}_{hk} - M (1-x_{hk})  &&  \forall h \in N_{d}, k \in N_{a}, k \neq h \label{eqn:secmodel_con21}\\
%------------------------------------%-----------------------------------------%
t^{d}_{j} &\geq t^{a}_{j} + s_{i} && \forall j \in C_{0} \label{eqn:secmodel_new1}\\
%------------------------------------%-----------------------------------------%
t^{d}_{j} &\geq t^{a}_{j} + f(b^{d1}_{j}, b^{a}_{j}) && \forall j \in S \label{eqn:secmodel_new2}\\
%------------------------------------%-----------------------------------------%
t^{'}_{j} &\geq t^{'}_{i}+c^{D}_{ij} +s_{j}-M \left(1- \sum_{\substack{k \in N_{a} \\\textlangle i,j,k \textrangle \in D} } y_{ijk}\right) && \forall j \in C^{'}, i \in N_{d}, i \neq j \label{eqn:secmodel_con22}\\
%------------------------------------%-----------------------------------------%
t^{'}_{k} &\geq t^{'}_{j}+c^{D}_{jk} -M \left(1- \sum_{\substack{i \in N_{d} \\ \textlangle i,j,k \textrangle \in D}} y_{ijk}\right) && \forall j \in C^{'}, k \in N_{a}, k \neq j \label{eqn:secmodel_con23}
\end{align}

Constraints (\ref{eqn:secmodel_con17})--(\ref{eqn:secmodel_con23}) are associated with the travel time of the two vehicles. In particular, constraints (\ref{eqn:secmodel_con17})--(\ref{eqn:secmodel_con20}) ensure that the travel time and drone range limit are correctly handled. Constraint (\ref{eqn:secmodel_con21}) indicates that if the EV travels from node $h$ to node $k$ where $h \in N_{d},k \in N_{a}$, its arrival time at node $k$ must incorporate its arrival time at node $h$, travel time from node $h$ to node $k$, the drone's launch time at node $h$ and retrieve time at node $k$. This constraint is not binding if the electric vehicle does not travel from node $h$ to node $k$. Constraints (\ref{eqn:secmodel_con22}) and (\ref{eqn:secmodel_con23}) are associated with the drone's arrival time. Suppose there is a drone route of $ \textlangle i,j,k \textrangle$, then the drone's arrival time of node $j$ and node $k$ should be related to the drone's travel time between $i$ and $j$,  $j$ and $k$ and the drone's retrieve time at node $k$. 

\par  %Besides, the retrieve time at node $k$ should also be considered.

\hfill \break
\noindent \textbf{Ordering Constraints:}
\begin{align}
t^{'}_{k} - t^{'}_{j} + c^{D}_{ij} & \leq Q_{d} + M (1- y_{ijk}) \notag \\
& \forall k \in N_{a}, j \in C^{'}, i \in N_{d}, \textlangle i,j,k \textrangle \in D  \label{eqn:secmodel_con24} \\
% *****************
u_{i} - u_{j} &\leq -1+(c+ms+2)(1-p_{ij}) && \forall i, j \in N^{'}, i \neq j \label{eqn:secmodel_con25}\\
% *****************
p_{ij} + p_{ji} &= 1 && \forall i, j \in N^{'}, i \neq j \label{eqn:secmodel_con26}\\
% *****************
t^{'}_{l} &\geq t^{'}_{k}- M \left(3-\sum_{\substack{j \in C^{'} \\ \textlangle i,j,k \textrangle \in D \\ j \neq l} }y_{ijk}- \sum_{\substack{m \in C^{'}\\m \neq i \\ m \neq k\\ m \neq l }}\sum_{\substack{n \in N_{a} \\ \textlangle l,m,n \textrangle \in D \\ n \neq i \\ n \neq k }}y_{lmn}-p_{il}\right) \notag \\
& \forall i, l \in N_{d}, k \in N_{a}, i \neq k \neq l, \textlangle i,j,k \textrangle, \textlangle l,m,n \textrangle \in D  \label{eqn:secmodel_con27}
\end{align}

Constraints (\ref{eqn:secmodel_con24})--(\ref{eqn:secmodel_con27}) are associated with ordering the two vehicles. Constraint (\ref{eqn:secmodel_con24}) ensures that the drone route should be within the drone's flight range. Constraint (\ref{eqn:secmodel_con25}) is a sub-tour elimination constraint and constraint (\ref{eqn:secmodel_con26}) ensures the correct ordering of two different nodes. Constraint (\ref{eqn:secmodel_con27}) indicates that if there exists two drone route deliveries $\textlangle i,j,k \textrangle$ and $\textlangle l,m,n \textrangle$ and node $i$ is visited before node $l$ by the EV, then node $l$ must be visited after node $k$. \par

\hfill \break
\noindent \textbf{Domain Constraints:}
\begin{align}
t_{0} &= 0   \label{eqn:secmodel_con28}\\
% *****************
p_{0j} &= 1 && \forall j \in N_{a} \label{eqn:secmodel_con29}\\
% *****************
x_{ij} &\in \{0, 1\}  && \forall i \in N_{d} \label{eqn:secmodel_con30}\\
% *****************
y_{ijk} &\in \{0, 1\}  && \forall i \in N_{d} \label{eqn:secmodel_con31}\\
% *****************
1 \leq u_{i} &\leq c+ms+2 && \forall i \in N_{a} \label{eqn:secmodel_con32}\\
% *****************
t_{i} &\geq 0 && \forall i \in N \label{eqn:secmodel_con33}\\
% *****************
t^{'}_{i} &\geq 0 && \forall i \in N \label{eqn:secmodel_con34}\\
% *****************
p_{ij} &\in \{0, 1\}  && \forall i \in N_{d}, j \in N_{a}, j \neq i \label{eqn:secmodel_con35}
\end{align}

Constraints (\ref{eqn:secmodel_con28})--(\ref{eqn:secmodel_con35}) specify the domain of all the decision variables.

\subsection{Modeling various partial recharge policies}
The EVTSPD-FF model proposed in \citep{zhu2022full} assumes that the EV is a Battery Electric Vehicle (BEV) and the charging station is a battery swapping station (BSS) so that the EV could refresh its battery energy to full capacity with a fixed time every time it visits the charging stations. However, in the real-world application, the full-charge assumption may not be applicable in specific scenarios. Below, we discuss some of the reasons why the partial-recharge assumption is preferable to the full-charge assumption.

Firstly, in a time-emergency routing scenario, e.g. a logistic company that aims to deliver parcels with time constraints, adopting a full-charge policy might render sub-optimal solutions. For example, when the EV visits a charging station near the end of its route, it does not need to charge its energy to full capacity to return to the depot. This situation also appears when the EV visits two consecutive charging station nodes, where full charging is unnecessary for the EV to travel to the next charging station. In this situation, more time-saving could be achieved by enabling partial recharge. \par
\begin{figure}[H]
    \centering
    \includegraphics[width = .6\textwidth]{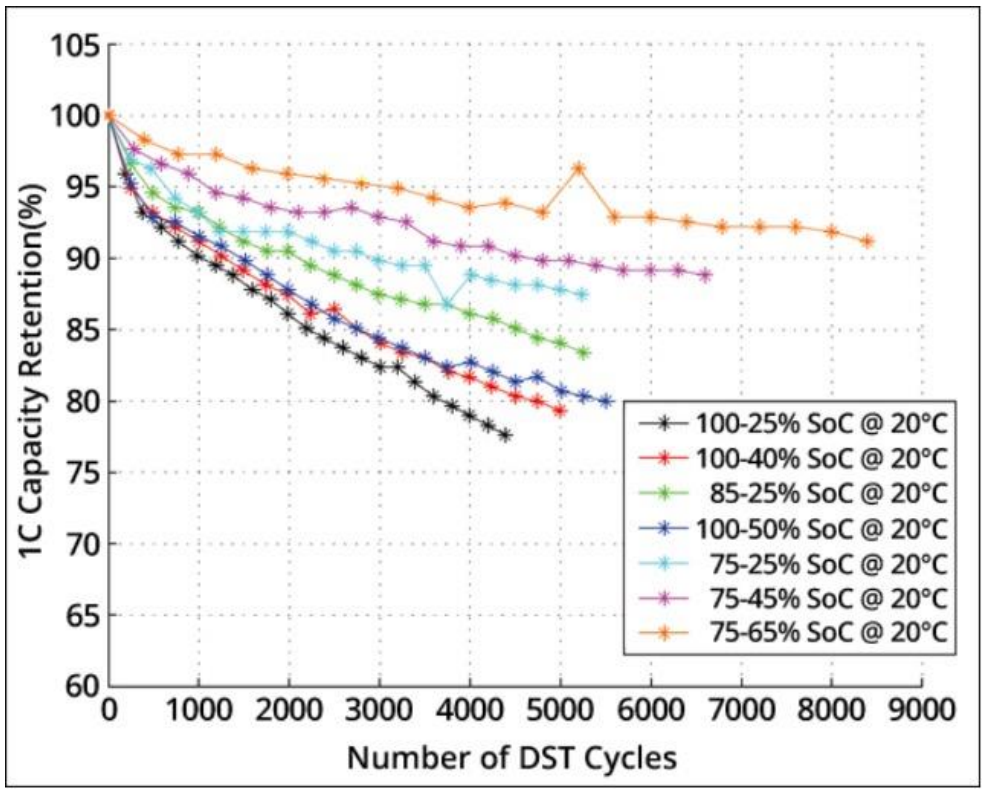}
    \caption{Battery capacity with respect to dynamic stress tests cycles (Source: \citep{BatteryUniversity})}
    \label{fig:DST_cycles}
\end{figure}

Secondly, in terms of vehicle maintenance, deep charge/discharge is detrimental to the useful life of lithium batteries which are commonly adopted in EVs and would significantly increase the battery's degradation speed, as indicated by \citep{BatteryUniversity}. Based on this research, a complete charge/discharge will fasten the battery's capacity fading process and lower the number of dynamic stress test cycles, as shown in Figure \ref{fig:DST_cycles}.

Last but not least, the number of available battery swap stations in the real world is fairly small compared to the number of available plug-in charging stations. This situation indicates that modeling all the charging station nodes in EVTSPD as BSSs might be impractical in a real-world application.\par

For these reasons, in this paper, we relax the fixed-time-full-charge assumption and instead assume that the EV is a plug-in hybrid electric vehicle. Meanwhile, the charging stations are connected to a plug-in electric grid, and the EV could be charged at CS nodes with any amount of energy. Both the amount of energy charged at CS nodes and the charging time, which depends on the initial SoC level and the amount of energy charged, are decision variables. This problem is named EVTSPD-P, where "P" stands for partial recharge. \par

\begin{figure}[H]
    \centering
    \includegraphics[width = .6\textwidth]{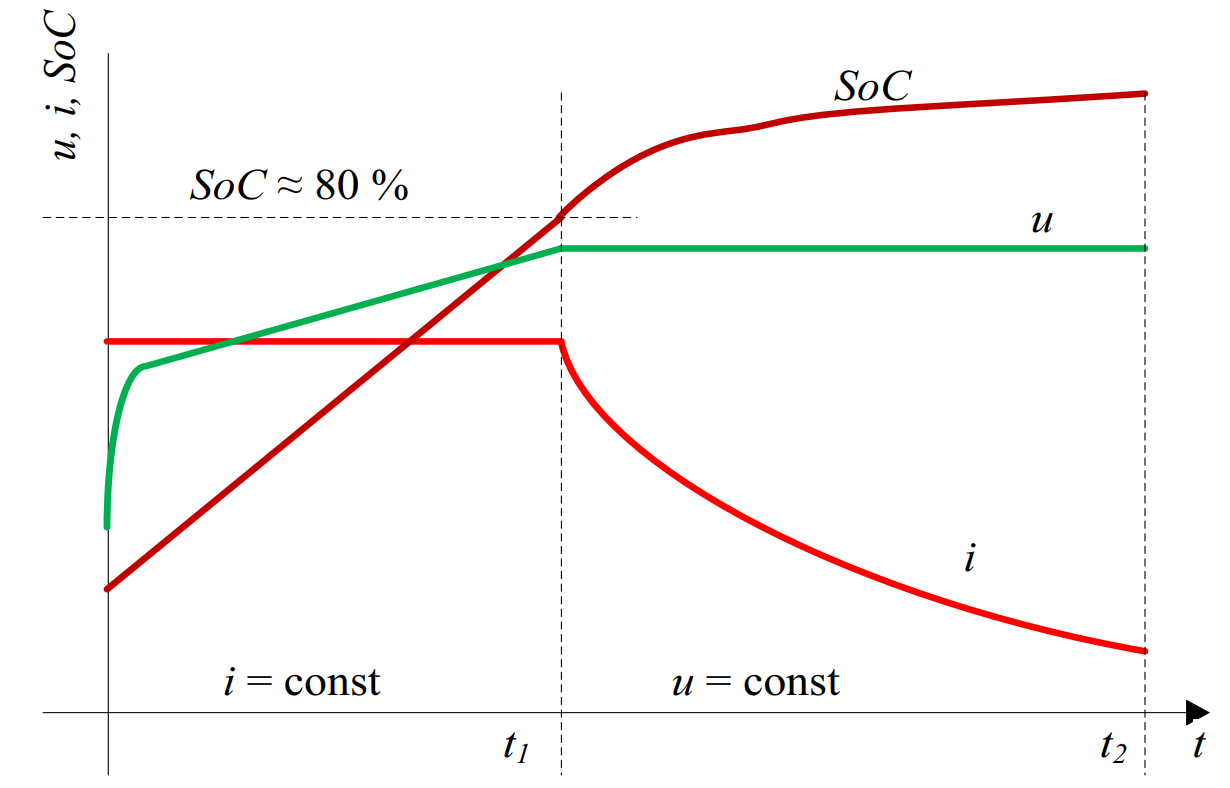}
    \caption{A typical charging curve (Source: \citep{hoimoja2012ultrafast}) }
    \label{fig:charging_curve}
\end{figure}

In the plug-in charging stations, the actual plug-in time-SoC curve is a non-linear concave function, as shown in Figure \ref{fig:charging_curve}. In the figure, $u$ and $i$ represent terminal voltage and current, respectively. The time-SoC curve is non-linear because the $i$ and $u$ change during the charging process. In the first stage, SoC increases linearly with charging time as the current $i$ is constant. After SoC reaches a specific value, the current $i$ decreases exponentially while the voltage $u$ is kept constant to avoid damage to the battery. As a result, the SoC grows concavely with time in the second phase. This charging curve indicates the charging efficiency decreases drastically when the SoC level is high, so it might not be worthwhile for an EV to charge a small amount of energy using a long period of time. 

Apparently, one key step in modeling EVTSPD-P is to estimate the concave time-SoC function. Multiple approximation methods have been used to construct an analytical expression to model this function, for example, \cite{wang2013estimation, bruglieri2014vehicle}, either by differential equation or linear approximation. \cite{zhang2020review} reviews the methods for estimating the time-SoC function. However, most of these techniques are too difficult to include in a MILP formulation or may cause deviations that cannot be ignored. Thus, we only consider the two most common approaches discussed below in this research.

\subsubsection{EVTSPD-P with linear time-SoC function estimation}
\begin{figure}[H]
    \centering
    \includegraphics[width = .6\textwidth]{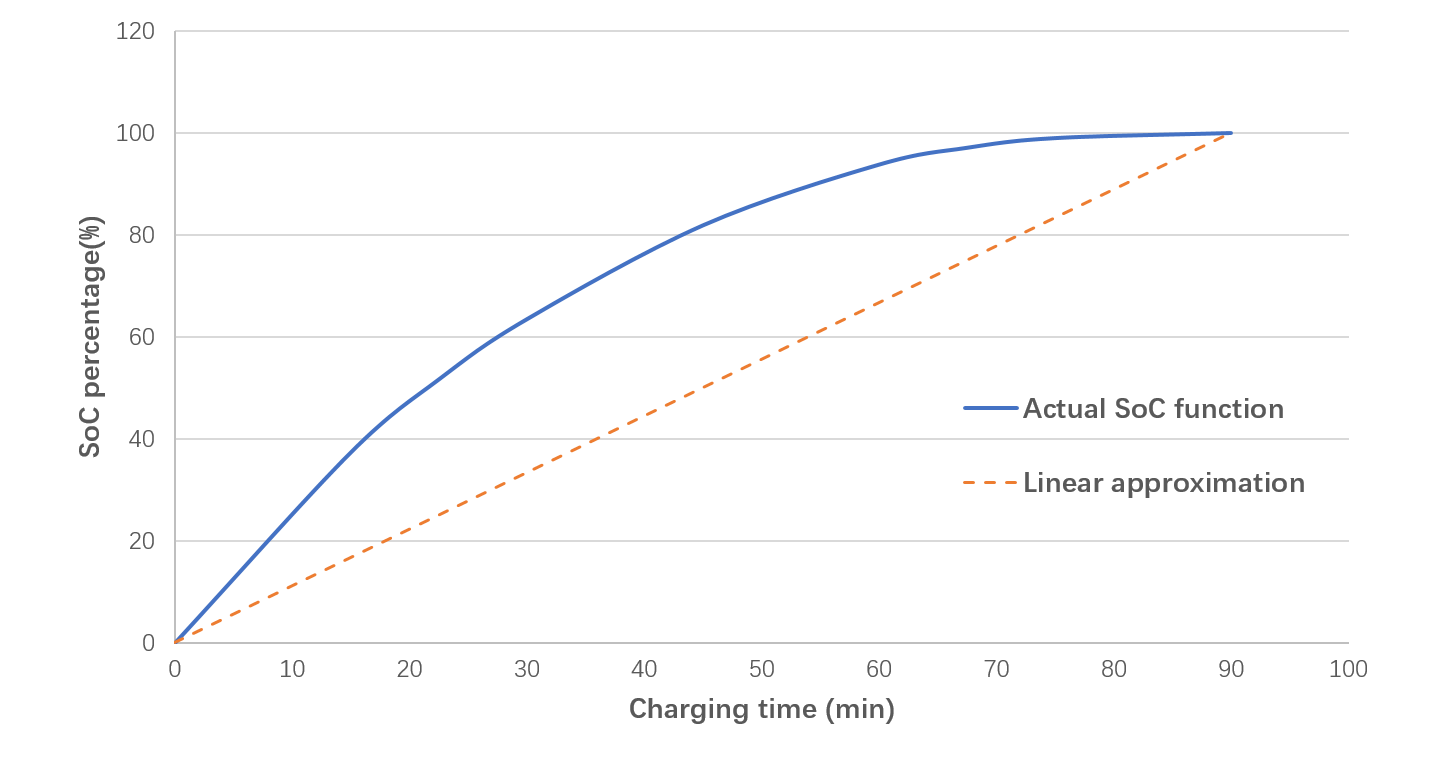}
    \caption{Actual time-SoC function and linear function approximation}
    \label{fig:SoC_linearApprox}
\end{figure}

In this subsection, we assume that the concave time-SoC function is approximated by a linear function,  $t(\Delta) = H * \Delta$ where $t(\Delta)$ is the time needed to charge $\Delta$ amount of energy and $H$ is a constant associated with each charging station. A representation of the time-SoC function and its linear approximation is shown in Figure \ref{fig:SoC_linearApprox}. Note that in this paper, we tend to use a conservative estimate of the actual time-SoC function, which indicates that for any charging time $t$, the actual SoC level after charging is higher than the linear approximation, to ensure the feasibility of the obtained solution. Similar assumption is made in \cite{felipe2014heuristic, sassi2014vehicle, bruglieri2015variable, schiffer2017electric, keskin2016partial}. In this paper, the EVTSPD model adopting linear SoC assumption is named EVTSPD-PL, where "P" stands for a partial charge, and "L" stands for linear approximation. \par

To incorporate the linear time-SoC function into the MILP model, the constraint (\ref{eqn:secmodel_new2}) is now represented as 
\begin{equation}
    t_{j}^{d} \geq t_{j}^{a} + H(b^{d1}_{j}-b^{a}_{j}), \forall j \in S \label{eqn:newcon6}
\end{equation}
and the rest of the formulation remains unchanged.

\subsubsection{EVTSPD-P with piecewise linear time-SoC function estimation}
\begin{figure}[H]
    \centering
    \includegraphics[width = .6\textwidth]{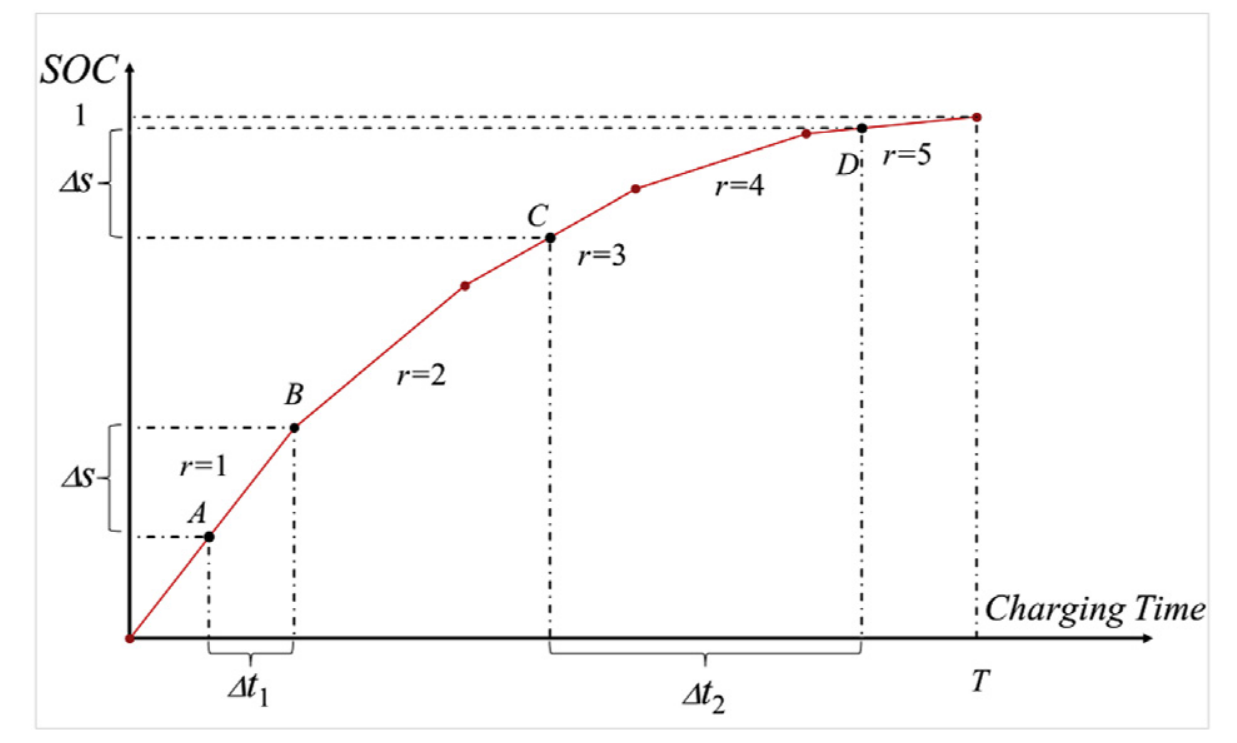}
    \caption{Relationship between charged time and charged energy amount}
    \label{fig:charge_amount}
\end{figure}

\begin{figure}[H]
    \centering
    \includegraphics[width = .6\textwidth]{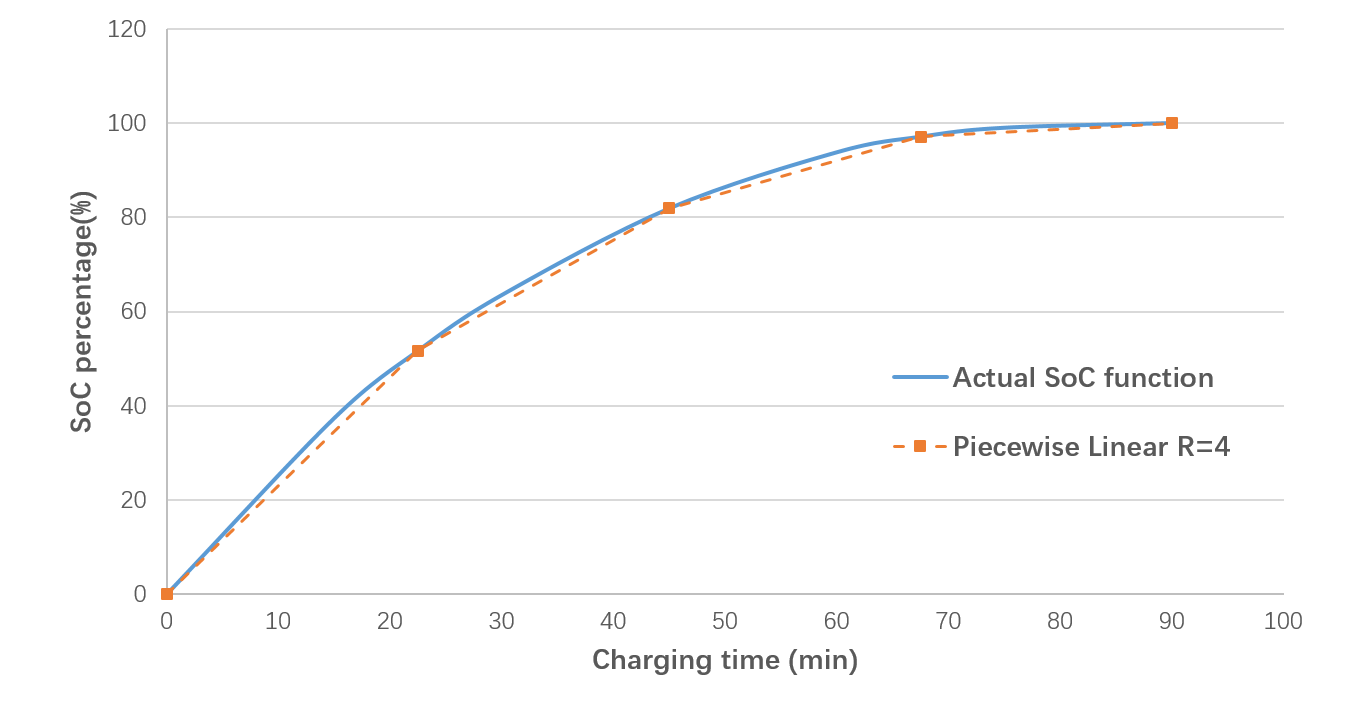}
    \caption{Actual time-SoC function and piecewise linear function approximation with $|R|=4$}
    \label{fig:SoC_piecewiseLinear}
\end{figure}

One major issue with the linear approximation of the time-SoC function is its accuracy. The relationship between the amount of charged energy and the needed charging time is not linear. It depends on the initial time-SoC level before charging, as represented in Figure \ref{fig:charge_amount}. \cite{Montoya2017} was the first to model the concave time-SoC function with a piecewise linear function. This is more precise than its linear counterpart, as represented in Figure \ref{fig:SoC_piecewiseLinear}, where $|R|$ represents the number of segments in the piecewise linear function. In \citep{Montoya2017}, four additional kinds of decision variables are needed in the MILP formulation. The total cost when using the non-linear charging function was found to be 2.70\% cheaper than that of the linear function on average. Later, this piecewise linear approximation is simplified by \cite{zuo2019new}, in which only one additional decision variable is needed. In this paper, we use the same technique as in \cite{zuo2019new} to model the piecewise linear time-SoC function. The EVTSPD model adopting the piecewise linear time-SoC function approximation is named EVTSPD-PP, where "PP" stands for "partial" and "piecewise linear".\par

Before we explain how to incorporate piecewise linear time-SoC function into the proposed three-index MILP model, some necessary changes in the notations are introduced below:
\begin{enumerate}
    \item $e^{T}_{ij} \in [0,1]$: energy cost relative to EV's energy capacity when traversing arc $(i,j)$.
    \item $e^{D}_{ijk}$: energy cost relative to EV's energy capacity when performing sortie $(i,j, k)$.
    \item $b^{a}_{i} \in [0,1]$:  remaining energy relative to EV's energy capacity upon arriving at node $i$.
    \item $b^{d1}_{i} \in [0,1]$:  remaining energy relative to EV's energy capacity upon departure at node $i$ while before launch operation.
    \item $b^{d2}_{i} \in [0,1]$:  remaining energy relative to EV's energy capacity upon departure at node $i$ after launch operation.
\end{enumerate}
In general, for the new notation, all the energy value associated with the node or arc are now represented as its ratio to the EVs battery energy capacity instead of its absolute value. Accordingly, the EVs battery capacity now equals one in the new notation.\par

Now we are ready to present the steps of modeling the piecewise linear time-SoC function. Denote $f(t)$ as the concave energy charging function of charging time $t$. It can be surrogated by a set of consecutive secant lines, as shown in Figure \ref{fig:SoC_piecewiseLinear}. Each secant line $r \in R$ is represented by slope $K_{r}$ and intercept $B_{r}$. Denote $S_{i}$ as energy charging level and $t_{i}$ as the required charging time if the battery is charged from $SOC=0\%$ to $SOC=S_{i}$ at charging station $i$. Thus, it is obvious that
\begin{equation}
    S_{i} \leq K_{r}t_{i}+ B_{r}, \forall r \in R \label{eqn:concave_inequality}
\end{equation}\par

Now, let $(t_{i}, s_{i}), i = 1,2,...,k+1$ be the breakpoints of the secant lines of SOC charging curve. Note that since there are a total of $|R|=k$ secant lines, there are $k+1$ breakpoints. Before charging, the current SOC level at charging station node $i \in S$ is $b^{a}_{i}$. Denote $t^{ba}_{i}$ as the time needed if the battery is charged from $SOC=0\%$ to $SOC=b^{a}_{i}$. Let $\Delta S$ be the amount of energy that needs to be charged and $\Delta t$ be the required charging time to charge the battery from $SOC =b^{a}_{i}$ to $SOC =b^{a}_{i}+\Delta S$. Then, with the same principle of Equation (\ref{eqn:concave_inequality}), we can derive the following equation:
\begin{equation}
    b^{a}_{i}+\Delta S \leq K_{r} *(t^{ba}_{i}+\Delta t) + B_{r}, \forall r \in R \label{eqn:actual_inequality1}
\end{equation}
Note that in our MILP formulation $\Delta t = t^{d}-t^{a}$ and $\Delta S = b^{d1}_{i} - b_{i}^{a}$. So the above equation is equivalent to 

\begin{equation}
    b^{d1}_{i} \leq K_{r} *(t^{ba}_{i}+t^{d}_{i}-t^{a}_{i}) + B_{r}, \forall r \in R, i \in S \label{eqn:actual_inequality2}
\end{equation}

Equation (\ref{eqn:actual_inequality2}) is what we need to model relationship between $b^{d1}_{i}$, $t^{d}$, and $t^{a}$. One problem here is that this equation has an extra variable $t^{ba}_{i}$. Thus, the following equations are needed to linearize equation (\ref{eqn:actual_inequality2}).\par

Denote $\alpha_{ri}$ as the binary variable that equals to one if a particular secant line $r$ is used to represent the relation between $b^{a}_{i}$ and $t^{ba}_{i}$. Thus, 
\begin{equation}
   b^{a}_{i} \geq s_{r} +(\alpha_{ri}-1), \forall r \in R, i \in S 
  \label{eqn:newcon1}
\end{equation}
\begin{equation}
   b^{a}_{i} \leq s_{r+1} +(1-\alpha_{ri}), \forall r \in R, i \in S
   \label{eqn:newcon2}
\end{equation}
\begin{equation}
    \sum_{r \in R}\alpha_{ri} = 1, \forall i \in S \label{eqn:newcon3}
\end{equation}
where constraints (\ref{eqn:newcon1}) and (\ref{eqn:newcon2}) ensure that if $\alpha_{r}=1$ then $b^{a}_{i}$ is bounded by $s_{r}$ and $s_{r+1}$ and constraints (\ref{eqn:newcon3}) ensures that only one secant line is selected where $SOC = b^{a}_{i} $ belongs to.\par
Now, we need to establish the relationship between existing variables with $t^{ba}$, which can be expressed as
\begin{equation}
    b^{a}_{i} \geq K_{r} t^{ba}_{i} +B_{r}-(1-\alpha_{r}), \forall r \in R, i \in S \label{eqn:newcon4}
\end{equation}
\begin{equation}
    b^{a}_{i} \leq K_{r} t^{ba}_{i} +B_{r} + (1-\alpha_{r}), \forall r \in R, i \in S \label{eqn:newcon5}
\end{equation}\par

With these new constraints, constraints (\ref{eqn:secmodel_new2}) now is replaced by
\begin{equation}
    t_{j}^{d} \geq t_{j}^{a} + s_{j}, \forall j \in S \label{eqn:newcon6}
\end{equation}
and the complete formulation for EVTSPD-PN is objective (\ref{eqn:secmodel_con1}) with constraints (\ref{eqn:secmodel_con2}) to (\ref{eqn:secmodel_new1}), constraints (\ref{eqn:secmodel_con22}) to (\ref{eqn:secmodel_con35}) and constraints (\ref{eqn:actual_inequality2}) to (\ref{eqn:newcon6}).

\subsection{Modeling additional side constraints in proposed MILP model}
This section shows how some of the commonly-encountered side constraints can be modeled based on the proposed three-index EVTSPD model. For example, when there are \textit{incompatible customers},  which can only be served by the EV and cannot be served by the drone, it can be considered by constructing only feasible drone sorties when creating drone sortie sets $D$. Similarly, when creating set $D$, we can also consider limited \textit{drone flying range} by only incorporating drone sorties whose length lies within the flying range. 

\subsubsection{Launch and retrieve times}
Define $S_{L}$ as the time needed to load and prepare the drone before the launch action and define $S_{R}$ as the time needed to unload and prepare the drone before the retrieve action. To model the launch and retrieve time, constraint (\ref{eqn:secmodel_con21}) is modified as constraint (\ref{eqn:secmodel_con35}) and constraint (\ref{eqn:secmodel_con23}) is modified as constraint (\ref{eqn:secmodel_con36}).
\begin{equation}
t^{a}_{k} \geq t^{d}_{h}+c^{T}_{hk} + S_{L}\sum_{\substack{l \in C^{'}}} \sum_{\substack{m \in N_{a} \\ \textlangle h,l,m \textrangle \in D} } y_{hlm} + S_{R}\sum_{\substack{i \in N_{d}}} \sum_{\substack{j \in C^{'} \\ \textlangle i,j,k \textrangle \in D} } y_{ijk} - M (1-x_{hk}),  \forall \ h \in N_{d}, k \neq h \label{eqn:secmodel_con35}\\
\end{equation}

\begin{equation}
t^{'}_{k} \geq t^{'}_{j}+s_{j}+c^{D}_{ij} + S_{R}-M \left(1- \sum_{\substack{i \in N_{d} \\ \textlangle i,j,k \textrangle \in D}} y_{ijk}\right), \ \forall j \in C^{'}, k \in N_{a}, k \neq j \label{eqn:secmodel_con36}
\end{equation}

As constraint (\ref{eqn:secmodel_con35}) allows $t_{j}{'}$ to be greater than the arrival time at customer $j$, our model corresponds to the ``wait'' model as explained in \cite{dell2019drone}, which indicates that the drone is allowed to wait at customer $j$'s location after serving $j$, while incurring no extra energy consumption in constraints (\ref{eqn:secmodel_con36}).

\subsubsection{Maximum number of customers per truck leg}
A \textit{truck leg} is defined as a concatenation of arcs traversed by the truck between the launch node and retrieve node of a drone sortie. One common side constraint is limiting the maximum number of customers within a truck leg. Denote $\Bar{n}$ as the maximum number of customers per truck leg. To account for this limitation, we need additional constraints (\ref{eqn:secmodel_con37}). 
\begin{equation}
y_{ijk} (u_{k}-u_{i}) \leq \Bar{n},  \ \ \forall \ \textlangle i,j,k \textrangle \in D  \label{eqn:secmodel_con37}
\end{equation}
As constraints (\ref{eqn:secmodel_con37}) is non-linear, the big-M method is adopted to modify it to the linear constraint. Define new decision variable $l_{i} = y_{ijk}u_{i}$ and $r_{k} = y_{ijk}u_{k}$. Constraints (\ref{eqn:secmodel_con37}) can be modified as follows, where $M = c+ms+2$ is an upper bound of $u_{k}, k \in N$:

\begin{align}
r_{k}-l_{i} &= \Bar{n}   \label{eqn:secmodel_con38}\\
% *****************
l_{i} &\leq M*y_{ijk} && \forall \ \textlangle i,j,k \textrangle \in D \label{eqn:secmodel_con39}\\
% *****************
l_{i} &\leq u_{i}  && \forall \ i \in N_{d} \label{eqn:secmodel_con40}\\
% *****************
l_{i} &\geq u_{i}-M*(1-y_{ijk}) && \forall \ \textlangle i,j,k \textrangle \in D \label{eqn:secmodel_con41}\\
% *****************
r_{k} &\leq M*y_{ijk} && \forall \ \textlangle i,j,k \textrangle \in D \label{eqn:secmodel_con42}\\
% *****************
r_{k} &\leq u_{k}  && \forall \ k \in N_{a} \label{eqn:secmodel_con43}\\
% *****************
r_{k} &\geq u_{k}-M*(1-y_{ijk}) && \forall \ \textlangle i,j,k \textrangle \in D \label{eqn:secmodel_con44}\\
% *****************
l_{i}, r_{k} & \geq 0  && \forall \ i \in N_{d}, k \in N_{a} \label{eqn:secmodel_con45}
\end{align}

\subsubsection{Weight-dependent drone flying range}
These side constraints assume that the drone's flying range is not a fixed value but a function of the weight of the parcels it carries when traversing an arc. In the proposed model, the drone sortie set $D$ is defined as the set containing all the feasible sorties. Under the weight-dependent flying range assumption, in a sortie $\textlangle i,j,k \textrangle$, the consumed energy when traversing arc $(i,j)$ is not an fixed value but a function $f_{(ij)}^{e}(w_{j})$ where $f_{(ij)}^{e}$ is the energy consumption function on arc $(ij)$ and $w_{j}$ is the weight of the parcel of customer $j$. Thus, we can redefine and recalculate feasible sortie set $D$ to account for the weight-dependent drone flying range assumption.

\subsubsection{Loops}
A \textit{loop} is defined as a feasible drone sortie whose launch node and retrieve node are the same customer node. During the loop $\textlangle i,j,i \textrangle$, the drone serves the customer $j$ while the truck waits for the drone at node $i$. An additional side constraint assuming that the loops are allowed in operation is also commonly seen in the literature.\par
Under the loop assumption, the truck could wait at the depot or a single customer node before the drone serves all the remaining customers, and then the truck returns to the depot. Accounting for this possibility is non-trivial in the proposed three-index model. This is because, in the original model, each node $i$ is assigned with a unique decision variable $u_{i}$, which represents the order of the node visited. To model self-loop for every customer node in the network, multiple dummy copies of the nodes in the network are needed making it practically intractable. This result is also shown in the numerical analysis section.

\subsection{A comparison to the alternative MILP model proposed in \citep{zhu2022full}}
This section briefly discusses the difference between the three-index MILP model proposed in this paper and the arc-based MILP model proposed in \citep{zhu2022full}. \par
First of all, these two models use two distinct approaches to address the same problem, which is the EV's driving range constraint. The three-index model addresses this problem by creating dummy copies of the charging station nodes so that each CS node can be visited multiple times, resulting in an augmented network with an increased number of nodes. On the other hand, the arc-based model addresses this problem by constructing a new multigraph where each arc represents a possible path that connects two customer nodes. This modeling approach results in a multigraph with an increased number of arcs, but fewer nodes, as CS nodes are no longer needed in the network. \par
Second, as two models adopt two different approaches when dealing with the driving range constraints, this leads to the same EVTSPD model with slightly different assumptions. For example, in the three-index model, the launch/retrieve node could be customer nodes or CS nodes, while in the arc-based model, this assumption does not stand as CS nodes are not included in the network. \par
Third, both models have their own merits and disadvantages when handling different additional side constraints. For example, in the three-index model, it is straightforward to incorporate additional serving sequence constraints where a customer $i$ should be served before customer $j$, which is fairly common in the pick-up and delivery application. This could be accomplished by adding additional constraints $u_{i} \leq u_{j}$ in three-index model. However, addressing this constraint is non-trivial in the arc-based model, where additional decision variables must be created. \par

\section{Adaptive large neighborhood search solution method}
Considering the enormous number of constraints and decision variables involved in EVTSPD-P, solving the problem of practical size using the MILP formulation via a commercial solver is inefficient. This section proposes a meta-heuristic solution method based on an adaptive large neighborhood search. The reasons for choosing this specific method are explained first, followed by a description of the parameter setting and iterative steps.\par

In the past, multiple researchers \citep{campuzano2021multi, DeFreitas2018} has shown that VNS could produce some of the best results when solving traveling salesman problem with drone. However, when solving the EVTSPD with the VNS method, an important issue needs to be considered, that is, maintaining the solution's feasibility during the search process. During the search process of VNS, the location of nodes and arcs are constantly exchanged to explore new solution spaces so that a local optimal could be obtained via various neighborhood structures. When VNS is applied to solve TSPD, it is easy to maintain the feasibility of the newly obtained solution generated with various neighborhood structures. Thus, many new feasible solutions can be generated in each iteration of VNS. \par

However, this is no longer the case in EVTSPD, where the extra constraint, the EV's driving range limit, must be respected. When applying the neighborhood structures to a feasible solution $s$, the resulting solution $s^{t}$, obtained by exchanging the location of nodes or arcs in $s$, might be infeasible due to the EV's driving range limit constraint. This property indicates that the feasibility of $s^{t}$ needs to be checked before the evaluation of its quality. This feasibility checking function is of linear computational effort. If the resulting solution $s^{t}$ is infeasible, an additional function that aims to restore the solution's feasibility, denoted as $g(s)$, needs to be called. Such function is of computational complexity $\mathcal{O}(NS)$ where $N$ is the number of nodes in the network, and $S$ is the set of charging station nodes even if only one charging station node is allowed to be inserted into the EV's route (which is not necessarily true for EVTSPD as, between every two customer visits, the EV is permitted to traverse multiple charging station nodes). Besides, if the feasibility of the current solution cannot be restored, the solution needs to be discarded, which indicates that the solution space is not fully explored and often leads to a final solution of poor quality. For this reason, when solving EVTSPD instances of large size or with tight EV driving range constraints, the final solution obtained via VNS is often of poor quality or unexpected long computational time. \par

For this reason, we propose an alternative meta-heuristic, large neighborhood search (LNS) \cite{shaw1998using, ropke2006adaptive} . Typically in LNS, there exist two different kinds of operators, destroy methods and repair methods. Destroy method aims to "destroy" the current solution by eliminating a specific target number of nodes, while the repair method tries to restore a feasible solution given the partial solution and eliminated nodes from the previous step. Both the destroy and repair method contains inherent randomness so that the search process can be diversified and a more feasible region could be reached. \par

Compared to VNS, the primary advantage of LNS in solving EVTSPD is that the number of iterations and newly found solutions are much smaller as it adopts more costly methods than VNS on obtaining new solutions. In this way, the number of times feasibility restoring function $g(s)$ is called decreases significantly, resulting in less overall computational complexity. This concept is suitable for solving problems with tight constraints where it is relatively difficult to generate a new feasible solution. In the past literature, \citep{sacramento2019adaptive} adopts the same method to obtain solutions of good quality. \citep{keskin2016partial} also adopt the same method in solving electric vehicle routing with partial recharge.\par

This research adopts some of the most commonly seen destroy methods and repair methods in the past literature and creates a new repair method specially designed for EVTSPD. This new repair method uses a "heavy" but effective constraint-programming (CP) based approach to maintaining the solution feasibility in the search process. This does not increase the computational complexity dramatically as the number of iterations in LNS is typically much smaller than that in VNS. Besides, we adopt an adaptive LNS method (ALNS) in our implementation instead of normal LNS. In the adaptive LNS, both the destroy and repair methods are chosen to modify the current solution based on their performance in the previous iteration, which is represented as a performance score and would be updated repeatedly during the search process. In ALNS, the method selection is based on the roulette wheel selection rule. Denote $A$ as the set contains all the available destroy methods and $B$ as the repair method set. At the start of iteration $j$, the performance score of method $i$ is denoted as $w_{ij}$, then the probability of method $i$ being chosen is 
\begin{equation}
p_{ij} = \frac{w_{ij}}{\sum_{i}{w_{ij}}}
\end{equation}
After each iteration, $w_{ij}$ would be updated based on their performance in the iteration as follows: 
\begin{equation}
w_{i,j+1} = \rho w_{ij} + \frac{\lambda}{\tau} (1-\rho) 
\end{equation}

In the formula, $\rho$ is the score shrinking rate, and $\lambda/\tau$ is the score for the method in iteration $j$, where $\lambda$ is an indication of the fitness of the newly found solution and $\tau$ is an indication of the compute time of the current iteration. Intuitively, the value of $\lambda/\tau$ is higher if we find a better solution with a lower cost or enter a search space that is not visited before in a short amount of time. The different situations of $\lambda$ and $\tau$ are shown in Table \ref{tab:alns_parameter}.

\begin{table}[H]
\caption{Situations regarding the performance of destroy and repair methods in a single iteration} 
\label{tab:alns_parameter}
\begin{center}
\begin{tabular}{m{1.2cm} m{12cm}} 
\toprule
Parameter & Description  \\
\hline
$\lambda_{1}$   & The method finds a global best feasible solution  \\
$\lambda_{2}$   & The method finds a feasible solution which is worse than the global best solution\\
$\lambda_{3}$  & The method fails to find a feasible solution  \\
$\tau_{1}$   & The current iteration's compute time is greater than one second \\
$\tau_{2}$   & The current iteration's compute time less greater than one second \\
\bottomrule
\end{tabular}
\end{center}
\end{table}

The degree of destruction, that is, how many nodes should be eliminated from the current solution, also affects the destroy method performance. In this research, a random number between 1 and the maximum destroy number $m_{s}$, whose calculation is explained later in this section, is generated as the target number of nodes that would be removed. \par

We control the acceptance probability of a new solution with a global time-varying parameter $T$ called temperature. This concept is frequently used in Simulated Annealing \citep{kirkpatrick1983optimization}. In this research, if a better solution is found in an iteration, the new solution is always accepted. If the new solution $s^{t}$ has a higher cost, the probability of it being accepted is 
\begin{equation}
p(x^{t}) = e^{\frac{f(s)-f(s^{t})}{T}}
\end{equation}
where $f(s)$ denote the objective value of solution $s$. \par

As the solution space is explored, the value of $T$ gradually decreases from the initial value to zero. Given a pre-specified time limit for the search phase, the value of $T$ is updated as 
\begin{equation} \label{eqn:temp_update}
T = T_{0}\left(1-\frac{t}{t_{m}}\right)
\end{equation}
where $t$ is the elapsed search time and $t_{m}$ is the pre-set time limit. To avoid the value of the denominator decreasing to zero in formula (\ref{eqn:temp_update}),  the search terminates when the value of $T$ is sufficiently small. In this research, the initial temperature is set to be 1000.
\par

The pseudo-code for the ALNS algorithm is given in Algorithm \ref{alg:alns}.

\begin{algorithm} 
  \caption{ALNS} \label{alg:alns}
  \hspace*{\algorithmicindent} \textbf{Input:} EVTSPD instance and preset parameters\\
  \hspace*{\algorithmicindent} \textbf{Output:} EVTSP-D solution
  \begin{algorithmic}[1]
      \State $s \gets InitialSolution()$
      \State $s^{\ast} \gets s$
      \State $noImpro \gets 0$
      \While{$t \leq t_{m}$}
          \State Choose a pair of destroy and repair method $(d(s), r(s))$ based on their performance score
          \State Get a temporary solution $s^{temp} \gets r(d(s))$ 
          \State Get new feasible solution $s^{t}$. If $s^{temp}$ is feasible, then  $s^{t} \gets s^{temp}$. Otherwise, $s^{t} \gets g(s^{temp})$ where $g(s)$ is the feasibility restoring function
          \State $T=T_{0}(1-t/t_{m})$
          \If {$random(0,1) < e^{\frac{f(s)-f(s^{t})}{T}}$}
              \State $s \gets s^{t}$
          \EndIf
          \If {$f(s^{t}) < f(s^{\ast})$}
              \State $s^{\ast} \gets s^{t}$
              \State $noImpro \gets 0$
          \Else
          \State $noImpro \gets noImpro+1$
            \If {$noImpro > noImproMax$}
              \State $s \gets s^{\ast}$
              \State $noImpro \gets 0 $
            \EndIf
          \EndIf
          \State Update performance score for destroy/repair methods
      \EndWhile \label{euclidendwhile}
      \State \textbf{return} $s^{\ast}$  
  \end{algorithmic}
\end{algorithm}

\subsection{Initial Solution}
This section presents the modified Clarke-Wright savings algorithm (MCWS), a greedy maximum saving algorithm proposed in \cite{Erdogan2012}, to construct an EV's route that satisfies the battery constraint. The basic steps of this heuristic are as follows:
\begin{enumerate}
    \item Create back-and-forth tours that start at the depot, visit a customer, and end at the depot.
    \item Check if these back-and-forth tours are feasible in terms of EV's driving limit and place all feasible tours into set $A$ and infeasible tours into the set $B$. 
    \item For each infeasible tour in the set $B$, calculate the cost of a charging station node insertion into the tour and adopt the least cost one that makes the tour feasible. Place the feasible tour into the set $A$ and discard all the infeasible ones.
    \item Compute the savings associated with each pair of feasible tours in the set $A$. The savings is defined as the decreased cost if two tours are merged. The savings are stored in a set $SPL$ in decreasing order.
    \item While $SPL$ is not empty, merge the tours in the $SPL$ and update the set until no merge is available. During the merge process, insert charging station nodes if necessary. Meanwhile, check if any charging station nodes could be removed from the tour, given this move does not result in the violation of battery constraint.
\end{enumerate}

Note that the resulting initial solution generated with MCWS has all the customers served by the EV and has no drone sorties. 

\subsection{Removal methods}
In each iteration, the removal operator aims to remove some customer nodes from the current solution $s$ and return a partial solution $s^{p}$ and a node set $P$, which stores all the nodes that need to be inserted back to $s^{p}$ later on. The current solution contains an EV's route $s_{t}$ and a set $s_{d}$ containing all the existing drone sorties. Both the nodes in $s_{t}$ and nodes served by one of the drone sorties in $s_{d}$ can be removed in this process.\par

Now we will briefly explain, given a current solution $s$, how the maximum number of nodes $m_{s}$ that could be removed from $s_{t}$, is calculated.

\begin{figure}[H]
    \centering
    \includegraphics[width = .4\textwidth]{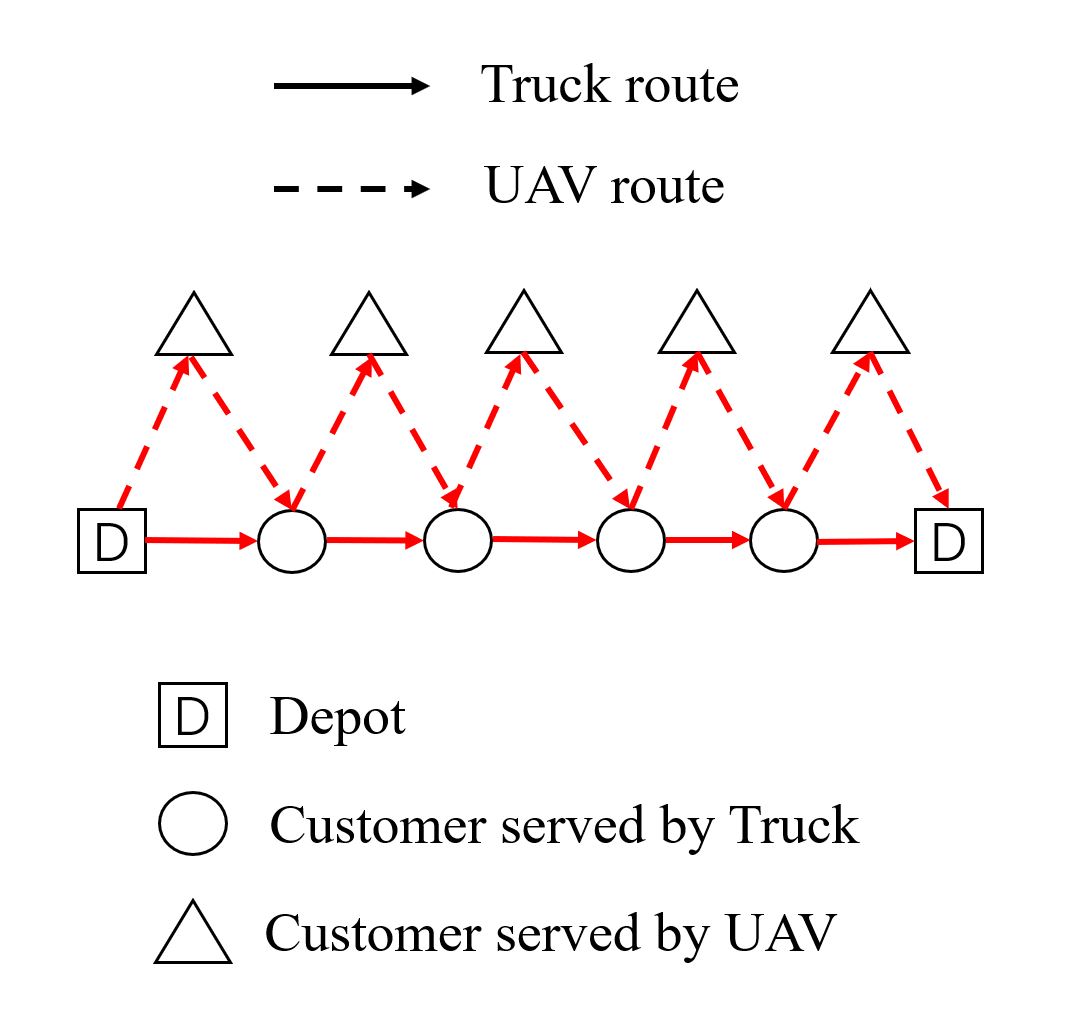} 
    \caption{Illustration of case when all drone sorties has no intermediate nodes}
    \label{fig:sortie_illu}
\end{figure}

Assume that the EV's route $s_{t}$ has a total of $n$ nodes (including the depot, customer nodes, and CS nodes) and that drone sortie set $s_{d}$ is empty. In this case, the maximum number of drone sorties that could be inserted into the EV's route is $n-1$. This value is also the number of traversed arcs in $s_{t}$. The maximum number of drone sorties is achieved when all the drone sorties have no intermediate nodes and are adjacent to each other, as shown in Figure \ref{fig:sortie_illu}. 

Assume for now that there are $d$ drone sorties stored in set $s_{d}$ and the maximum number of nodes $m_{s}$ are to be removed from $s_{t}$. If, later on all the $m_{s}$ nodes are served with the drone in the repairing process, the number of drone sorties in the resulting solution is $d+x$. This number should no greater than  $n-m_{s}-1$, which is the maximum number of drone sorties that could be inserted into $s_{t}$ after $m_{s}$ nodes are deleted from it. Setting $n-m_{s}-1$ equals to $d+m_{s}$ indicates that the corresponding $x=(n-1-d)/2$. In this paper, a randomly generated integer $\beta$ that lies between 1 and $\lfloor{(n-1-d)/2}\rfloor$ is regarded as the target number of customer nodes that would be deleted from the current solution. This approach aims to guarantee all the removed nodes can be served by the drone in the repairing process.\par

In this research, two different destroy methods are adopted. In every iteration, one of the destroy methods is chosen based on their performance score to modify the current solution $s$. These two destroy methods are explained next.

\subsubsection{Random removal}
The random removal method randomly removes $\beta$ customer nodes from the current solution. This method randomly selects and deletes a target number of customer nodes from the current solution $s$ and returns a partial solution and node-set. In the particular case when the selected node in $s_{t}$ also serves as a launch or retrieve node in one of the drone sorties in $s_{d}$, the corresponding drone node in that sortie will also be removed from the current solution.

\subsubsection{Cluster removal}
In the cluster removal method, a random customer node is initially chosen and removed from the current solution. We then progressively remove the remaining customer node closest to the previously chosen node until a target number of customers has been selected. Similarly, if the chosen node is a launch or retrieve node, then the corresponding drone node will also be removed.

\subsection{Repair Methods}
In each iteration of ALNS, with given partial solution $s^{p}$ and  node-set $P$, the repair method serves all the nodes in $P$ by inserting them into $s^{p}$ and returns a new feasible solution $s^{t}$, whose quality is hopefully better than the current solution $s$. This process needs to be carefully implemented to ensure that the resulting new solution $s^{t}$ is feasible. To accomplish this goal, several different heuristics are implemented. In this research, the repair method is a two-phase modification process. In general, the repair method first tries to construct a feasible solution by inserting picked nodes (and charging station nodes if necessary) into the current partial solution to achieve a feasible route. Then, different improvement methods are used in the second phase to improve the feasible solution by utilizing the drone to achieve more time savings. The main difference between different repair methods lies in the second improvement phase. As this method construct a feasible EV's route first, it can also be regarded as a variant of truck-first-sortie-second method. \par

In this paper, we consider three different repair methods, two of which are inspired from \cite{sacramento2019adaptive} and the last one is the newly proposed method. The three repair methods are explained below.

\subsubsection{Greedy repair}
The greedy repair method is the "cheapest'' of all three repair methods. This process can be divided into two phases. In the first phase, several feasible routes are constructed by inserting all the nodes in the picked node-set $P$ into the available truck's routes with a greedy method. All nodes are inserted at a location that yields the least extra cost. Later, this method tries to improve the solution by deploying the drone as much as possible.

In the second phase, this method tries to improve the current solution by checking every truck's node to see if the drone could be served if it is currently available. When creating the drone sorties, the previous node and the next node of the current node are regarded as the ``launch node'' and ``retrieve node,'' respectively. This way, the newly created sorties have no intermediate nodes, and the drone could serve customers as much as possible. If the resulting drone sortie satisfies the drone's flight limit constraint, the sortie is accepted, and we move to the next node in the truck's route. This process is repeated until every node in the truck's route has been checked. The detailed procedure is shown in Algorithm 2. 

\begin{algorithm}[H]
  \caption{Greedy repair}\label{euclid}
  \hspace*{\algorithmicindent} \textbf{Input:} partial solution $s^{p}$ and picked node set $P$\\
  \hspace*{\algorithmicindent} \textbf{Output:} an improved feasible solution of $s^{p}$
  \begin{algorithmic}[1]
      \For{node $i$ in $P$}
        \State Find the insert location with least extra cost among all the insertion locations;
        \State If the resulting coordinated route is feasible, adopt the move. Otherwise, insert CS to restore its feasibility.
      \EndFor
        \For{node $i$ in $r^{T}$}
          \If{$i$ could be served by drone with previous node as launch node and next node as retrieve node} 
            \State Add sortie $s$ to the current route;
            \State Delete $i$ from the truck's route.
          \EndIf
        \EndFor
      \State \textbf{return} $s^{p}$
  \end{algorithmic}
\end{algorithm}

\subsubsection{Nearby repair}
The second repair method, nearby repair, requires more effort than the greedy one. The first phase of the nearby repair is similar to that of the greedy repair, which is constructing a feasible truck's route by inserting the picked nodes in $P$ to the truck's route while ensuring the resulting truck's route is feasible. In the second improvement phase, the truck's route is first divided into several segments based on the current drone sortie. Each route segment is defined as "drone feasible" if the drone is on board and available to use in this segment. The truck node in such a feasible segment is called a "feasible" truck node. The goal of the improvement phase is to serve these feasible truck nodes smartly to yield more time-saving. Unlike the greedy method, which aims to use a drone to serve as many customers as possible, in the nearby repair, we aim to explore all possible drone sorties within the segment and pick the sortie with the highest positive time-saving. This process is repeated until no feasible segment is left or there is no potential drone sortie with positive time-saving. The computational complexity of the improvement operator is $o(|N|^{2})$ where $|N|$ is the number of nodes in the truck's route. The detailed process of the segment repair is shown in Algorithm \ref{alg:segment_repair}. 

\begin{algorithm}[H]
  \caption{Nearby repair}\label{alg:segment_repair}
  \hspace*{\algorithmicindent} \textbf{Input:} partial solution $s^{p}$ and picked node set $P$\\
  \hspace*{\algorithmicindent} \textbf{Output:} an improved feasible solution of $s^{p}$
  \begin{algorithmic}[1]
      \For{node $i$ in $P$}
        \State Insert node into $s^{p}$ before a customer node that is closest to the chosen node
      \EndFor
        \For{each feasible segment $f$ in $r^{T}$}
           \For{node $i$ in $f$}
             \State Examine every feasible drone sortie that serves $i$;
             \State Select sortie $s$ with the greatest saving or with lowest traverse node;
             \State If the resulting route is feasible, add sortie $s$ into $s^{p}$;
             \State Break $f$ into two remaining smaller segment and add them into segment set.
           \EndFor    
        \EndFor
     \State \textbf{return} $s^{p}$
  \end{algorithmic}
\end{algorithm}

\subsubsection{CP repair}
The third repair method, constraint programming (CP) repair, is the most computationally demanding of all three methods. In the first two methods, we aim to serve the truck nodes by identifying new drone sorties which yield high time-saving, while the current drone sorties in the partial solution $s^{p}$ remain unchanged. In the CP method, instead of keeping all the existing drone sorties in the partial solution unchanged, we will eliminate all these existing drone sorties, along with some extra customer nodes from the truck's route, and re-insert these customer nodes into the truck's route of the partial solution, by solving a scheduling problem. CP method is "heavier" than the two repair methods mentioned above in that it destroys all the remaining drone sorties and the truck routes in the partial solution. Besides, while the first two repair methods insert the customer nodes with the greedy method, in the CP repair method, we aim to construct new drone sorties optimally by solving a scheduling problem. Thus, CP method is more likely to escape local optima in the searching process, of course, with the downside that it requires more computational effort.

In this research, after CP method removes the existing drone sorties and several customer nodes in the truck's route, the second improvement phase is formulated as a scheduling problem in CP and solved by the CPLEX CP optimizer. The detailed process of the CP repair method is shown in Algorithm \ref{alg:CP_repair}. In line 8 of the algorithm, we aim to re-construct the solution using the selected nodes and partial solution. This subproblem is first proposed in \cite{yurek2018decomposition} and a detailed description of the MILP formulation is shown in the Appendix. The CP counterpart of this subproblem is described below.

\begin{algorithm}[H]
  \caption{CP repair}\label{alg:CP_repair}
  \hspace*{\algorithmicindent} \textbf{Input:} partial solution $s^{p}$ and picked node set $P$\\
  \hspace*{\algorithmicindent} \textbf{Output:} an improved feasible solution of $s^{p}$
  \begin{algorithmic}[1]
      \For{node $i$ in $P$}
        \State insert node into $s^{p}$ with greedy method.
      \EndFor
     \State $p' \gets$ null set;
     \State delete all drone sorties in $r_{i}^{D}$, add drone nodes into $p'$;
     \State randomly select truck nodes in $r_{i}^{T}$ and add them into $p'$;
     \State $r^{T'}_{i} \gets$ remaining truck route;
     \State solve the CP model or MIP model with $r^{T'}_{i}$ and $p'$.
      \State \textbf{return} $s^{p}$
  \end{algorithmic}
\end{algorithm}

\vspace{0.5cm}
\textbf{Indexes, sets and parameters used in the subproblem:}
{\renewcommand\arraystretch{1.0}
\noindent\begin{longtable*}
{@{}l @{\quad:\quad} p{15cm}@{}}
% {@{}>{\raggedright}p{2cm}lp{7.3cm}@{}}
$C$ & set of all customer nodes\\
$D$ & set of customers that need to be served by the drone \\
$n_{D}$  & number of customer nodes that needs to be inserted into the truck's route \\
$n_{T}$  & number of nodes in truck's route\\
$r_{T}$  & existing truck's route \\
$i$  & index of customer nodes $i \in C$\\
$j$ & index of insertion type $(l,r)$, where $l,r$ are the index of the launch/retrieve node in the truck's route with the constraint that $l < r \leq n_{T}$.   \\
$w_{ij}$  & incurred waiting time of inserting node $i$ by insertion type $j$ \\
\end{longtable*}}

\textbf{Decision Variables:}
{\renewcommand\arraystretch{1.0}
\noindent\begin{longtable*}
{@{}l @{\quad:\quad} p{15cm}@{}}
% {@{}>{\raggedright}p{2cm}lp{7.3cm}@{}}
$visit_{i}$ & interval variable representing the insertion of $i$-th customer node.\\
$visit_{ij}^{Alt}$ & optional interval variable representing the insertion of $i$-th customer node by type $j$. \\
$seq$ & sequence variable representing the order of all the existing $visit_{ij}^{Alt}$. \\
\end{longtable*}}

The subproblem can be modeled in the form of CP as follows:

\begin{itemize}
    \item \textbf{Objective}: Minimize the extra waiting time introduced by the synchronization of adding drone sortie to the existing truck's route.
    \item \textbf{C1}: Every customer nodes that is not is the truck's route should be served by the drone for exactly one time.
    \item \textbf{C2}: The drone route should within drone's driving range.
    \item \textbf{C3}: The drone should be launched and retrieved at the different customer nodes.
    \item \textbf{C4}: There is only one drone on the truck, so if the drone is launched, then before the truck retrieves the drone, it cannot be launched again.
\end{itemize}

Now, the objective of minimizing the total waiting time can be defined as follows:\par
\begin{equation}
Minimize \sum_{i} \sum_{j} \{\textit{\textbf{presenceOf}} (visit_{ij}^{Alt}) * w_{ij}\} 
\end{equation}
Since $ w_{ij}$ represents the waiting time of inserting node $i$ by insertion type $j$ and \textit{\textbf{presenceOf}} $(visit_{ij}^{Alt})$ is the binary indicator of whether $visit_{ij}^{Alt}$ exists, the $\sum_{i} \sum_{j} \{\textit{\textbf{presenceOf}} (visit_{ij}^{Alt}) * w_{ij}\}$ is the total waiting time, which is the objective function the CP aims to minimize.

Formulating C1, C3, and C4 can be seen as follows:
\begin{equation}
\textit{\textbf{Alternative}} (visit_{i}, visit_{ij}^{Alt}) \  \forall i \in D 
\end{equation}

\begin{equation}
\textit{\textbf{noOverlap}}(seq)
\end{equation}

$\textit{\textbf{Alternative}} (visit_{i}, visit_{ij}^{Alt})$ is a constraint function stating that for a customer node $i$, if $visit_{i}$ exists, then only one of the $visit_{ij}$ for all $j$ can exist. This constraint indicates that for every drone customer node, only one of the insertion type $j$ should exist. $\textit{\textbf{noOverlap}}(seq)$ states that the interval variables in the sequence variable $seq$ should not overlap, which corresponds to constraint C4, indicating that the drone's route should not overlap each other. Additionally, constraint C2 could be implemented in the period of data processing by eliminating all the $(i,j)$ combinations that exceed the drone's flight endurance when creating the decision variable $visit^{Alt}_{ij}$. 

\subsection{Feasibility resorting function $g(s)$ and cost evaluation function $f(s)$}
At each iteration, after the removal methods and repair methods are applied to the current solution $s$, we can get a temporary (possibly infeasible) solution $s^{temp}$. If $s^{temp}$ is feasible, then we have achieved a new feasible solution whose objective value is ready to be calculated. Otherwise, the feasibility resorting function $g(s)$ needs to be called and applied to $s^{temp}$ to try return a feasible solution. Note that function $g(s)$ may fail to obtain a new feasible solution. In this case, the temporary solution $s$ is discarded, and ALNS moves to next iteration. \par

\begin{figure}[H]
    \centering
    \includegraphics[width = .8\textwidth]{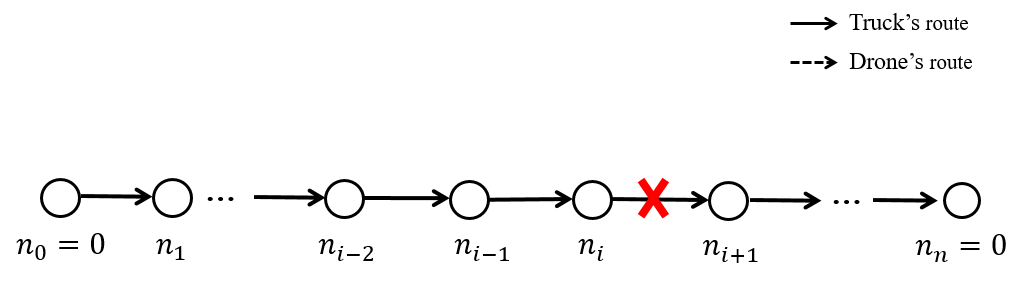}
    \caption{Comparison between ALNS and VNS }
    \label{fig:insertCSfunc_illu}
\end{figure}

As explained earlier, function $g(s)$ aims to insert charging station nodes into the EV's route of $s^{temp}$ to restore its feasibility. (Doing so would not affect the feasibility of drone sorties in $s^{temp}$ as we assume the drone can wait for the EV at the retrieve node.) Denote $s^{temp}_{t}$ as the EV's route in $s^{temp}$. Assume for now that the EV's route is as shown in \ref{fig:insertCSfunc_illu} where the EV's driving limit constraint is violated when the EV traverses arc $(n_{i}, n_{i+1})$. In this case, the detailed steps of feasibility resorting function $g(s)$ are as follows:
\begin{itemize}
    \item \textbf{Step 1}: identify the violating arc $(n_{i}, n_{i+1})$. If none found, the solution $s$ is feasible. Otherwise, move to Step 2.
    \item \textbf{Step 2}: set current arc $a^{cur}$  as arc $(n_{i}, n_{i+1})$. Identify the CS node $s_{i}$ that is closest to node $n_{i}$ and insert it into $a^{cur}$. If after inserting $s_{i}$ the EV could traverse the current arc successfully, return to Step 1. Otherwise, move to Step 3.
    \item \textbf{Step 3}: set current arc $a^{cur}$ as arc $(n_{i-1}, n_{i})$, which is prior to the current arc. Identify the CS node $s_{i-1}$ that is closest to node $n_{i-1}$. Try to insert $s_{i-1}$ OR $s_{i}$ into the current arc. If after inserting $s_{i}$ the EV could traverse the current arc successfully, return to Step 1. Otherwise, repeat Step 3, until we reach a preset limit.
\end{itemize}

Note that the main difference between Step 2 and Step 3 is that in Step 3 we try to insert node $s_{i-1}$ OR $s_{i}$ into arc $(n_{i-1}, n_{i})$ where in Step 2 we only try to insert node $s_{i}$ into arc $(n_{i}, n_{i+1})$. \par

If in one iteration of ALNS, we have obtained a feasible solution $s^{t}$, the quality or the objective value of this solution needs to be evaluated. The objective value of $s^{t}$ is simply the route completion time. As in EVTSPD-P, the EV is partially recharged in the charging station nodes, the charging time at each visited CS station needs to be calculated first. This calculation process is straightforward in EVTSPD-P. We can simply consider the extreme case, where each time the EV arrives at a CS node or the depot, its remaining energy is zero. In this way, the minimum energy required for an EV to travel from CS node or the depot to the next CS node or the depot in the route can be calculated. Back-tracking all the CS nodes in the EV's route enables us to get the minimum starting energy of the EV at the starting depot. As we assume that the EV starts its trip with full-capacity energy, we can add a fixed amount of energy to each CS node in the route if the minimum starting energy is less than the EV's energy capacity and then calculate the charging time at each CS nodes. Note that this approach guarantees the minimal charging time of the solution because the time-SoC function is linear or concave. \par

\section{Computational Experiments}
This section first compares the ALNS method proposed in this paper with VNS. Then the computational efficiency and accuracy of ALNS are tested by comparing it with solving the three-index MILP formulation with a commercial solver. A sensitivity analysis of the number of line segments used in the time-SoC function approximation process is also performed on randomly generated instances and real-world networks.\par

\subsection{Experimental setting}
In this paper, the numerical analysis is conducted on randomly generated instances with realistic parameter settings. For all the generated instances, the depot is located at $(0, 0)$ and the coordinates of the customer nodes and the charging station nodes are randomly generated between -20 km and 20 km. The distance matrix for the EV is calculated using the Manhattan metric, while the Euclidean metric is used for the UAV. Whenever each instance is generated, a simple MCWS algorithm firstly proposed in \citep{Erdogan2012} is used to solve the instance as EVTSP to guarantee that there exists a feasible solution to EVTSPD-P.\par

Similar to \citep{zhu2022full}, we assume the EV in the problem is Alke model ATX340E and the drone is DJI MATRICE 600 PRO. The parameter setting associated with the EV and drone is mostly adopted from the company's official website (\citep{Alke, DJI}), which is shown in Table \ref{tab:test_para}.

\begin{table}[H]
\begin{center}
\caption{The parameters Associated with the EV and the drone} \label{tab:test_para}
\begin{tabular}{ m{10cm} m{3cm} } 
\toprule
\textbf{Parameter}  & \textbf{Value} \\
\hline
\textbf{For EV (Alke TX340E with 10kWh battery):} & \\
Driving speed & 40 km/h \\
Maximum driving range (time) & 100 km (2.5 h) \\
Charging time at CS node & 2 min \\
\\
\textbf{For drone (DJI MATRICE 600 PRO):} & \\
Flying speed & 60 km/h \\
Maximum flight range (time) & 20 km (20 minutes) \\
Drone-EV energy consumption rate ratio & 0.4 \\
Payload capacity & 6 kg \\
Launch time & 100 s \\
Retrieve time & 20 s \\
\bottomrule
\end{tabular}
\end{center}
\end{table}

The proposed MILP formulation is implemented in Pyomo and solved with ILOG's CPLEX Concert Technology solver (version 12.6.3), while the ALNS algorithm is coded in python. For the CP repair method, the improvement subproblem is solved with the CPLEX CP optimizer (version 12.6.3). All experiments are run on a 3.6 GHz Intel Core i7 desktop with 32 GB RAM. For the test results, all the computational times are measured in seconds.

\subsection{Estimation of the concave time-SoC function}
\begin{figure*}
\centering
\begin{subfigure}[b]{0.475\textwidth}
    \centering
    \includegraphics[width=\textwidth]{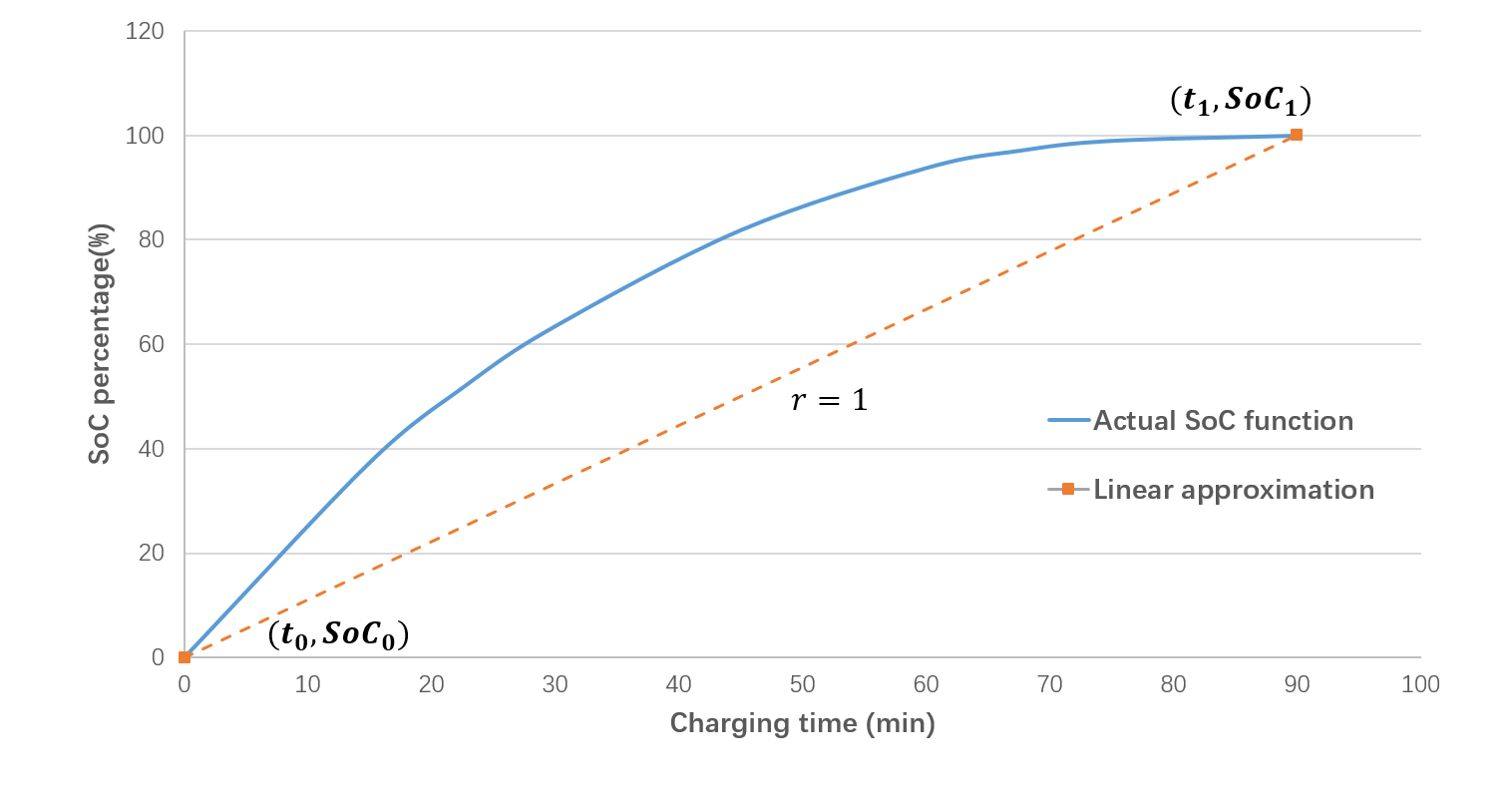}
    \caption[]%
    {{\small Linear approximation}}    
    \label{fig:SoC_r1}
\end{subfigure}
\hfill
\begin{subfigure}[b]{0.475\textwidth}  
    \centering 
    \includegraphics[width=\textwidth]{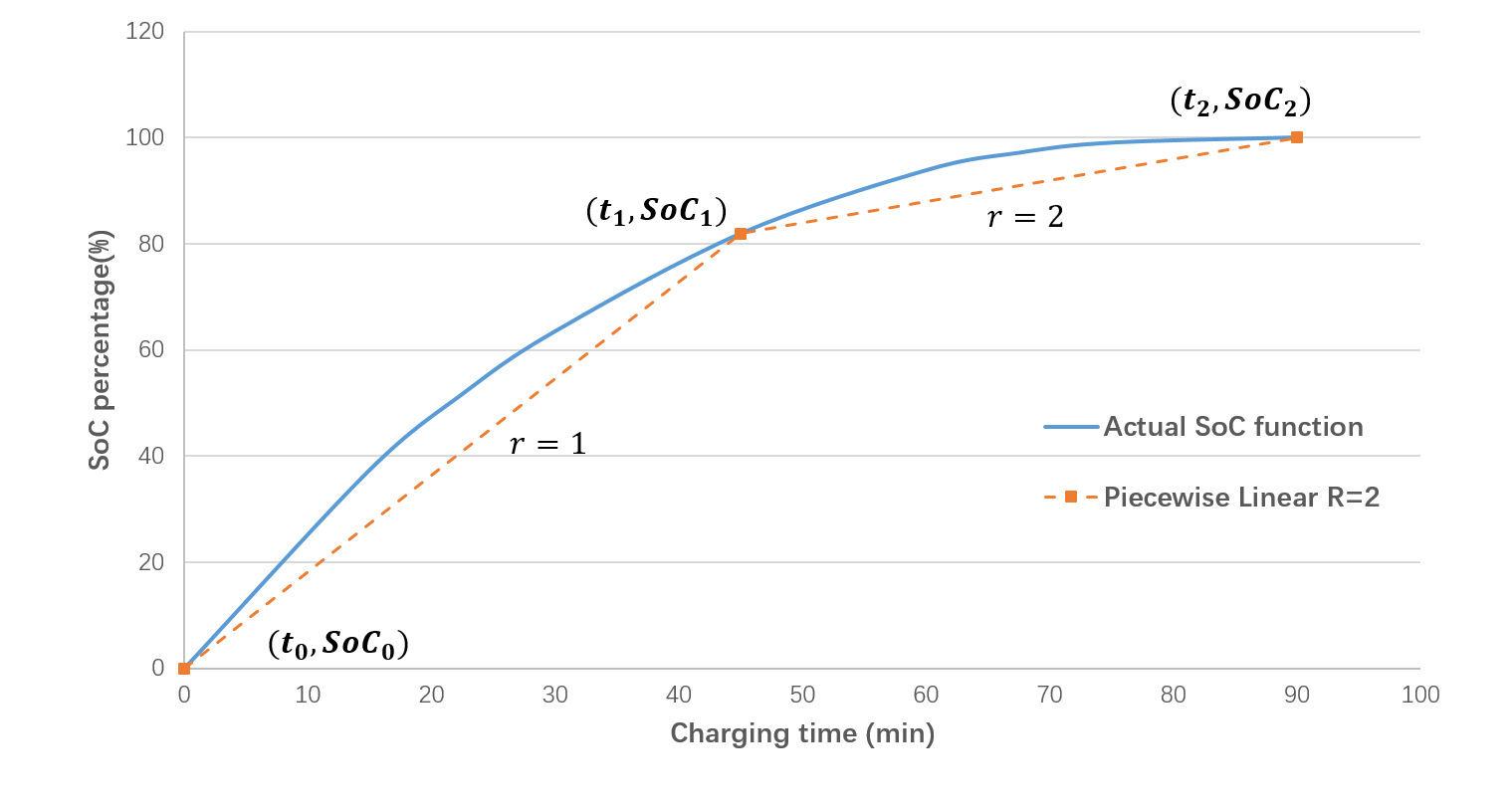}
    \caption[]%
    {{\small Piecewise linear approx R=2}}    
    \label{fig:SoC_r2}
\end{subfigure}
\vskip\baselineskip
\begin{subfigure}[b]{0.475\textwidth}   
    \centering 
    \includegraphics[width=\textwidth]{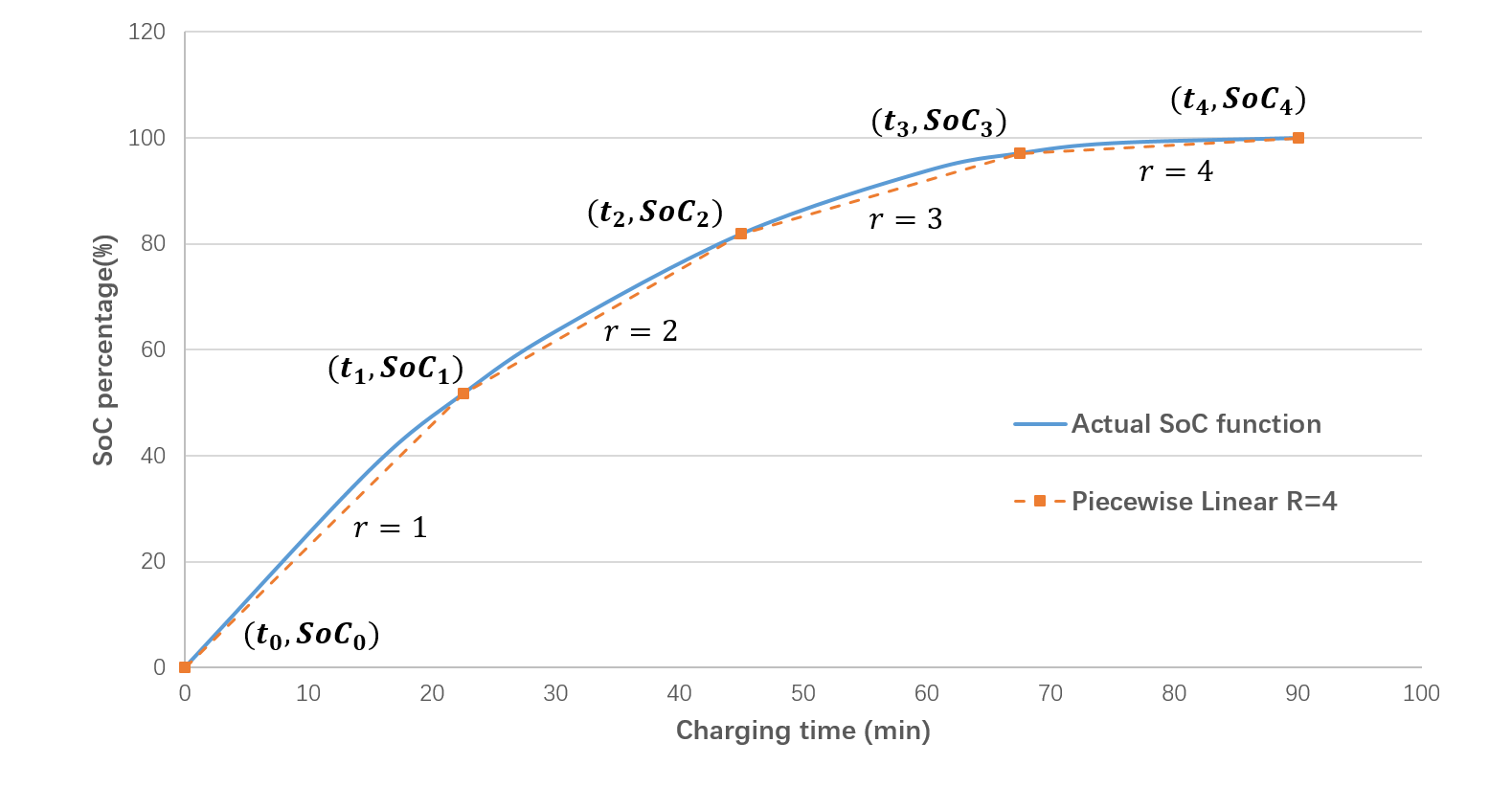}
    \caption[]%
    {{\small Piecewise linear approx R=4}}    
    \label{fig:SoC_r4}
\end{subfigure}
\hfill
\begin{subfigure}[b]{0.475\textwidth}   
    \centering 
    \includegraphics[width=\textwidth]{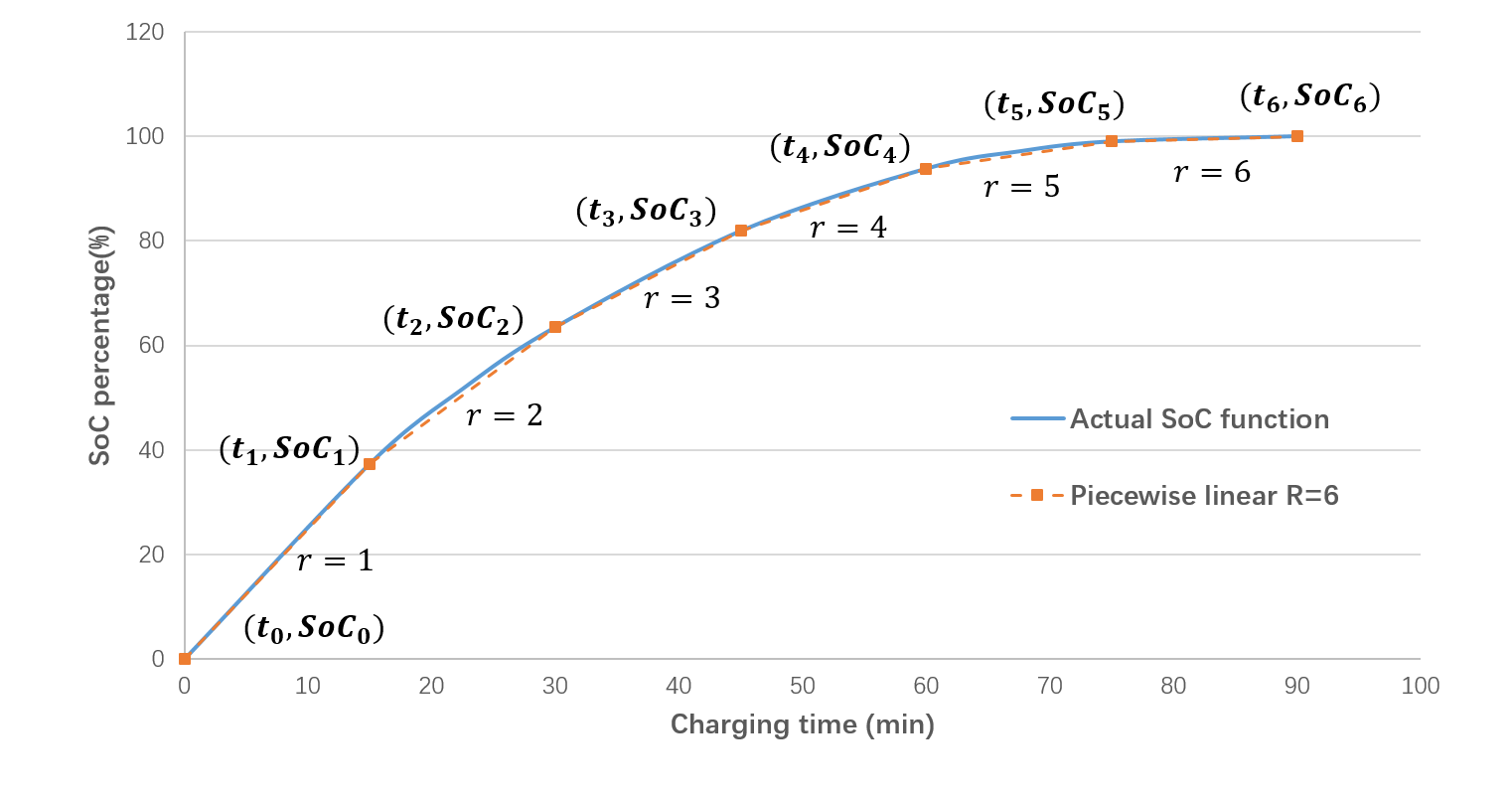}
    \caption[]%
    {{\small Piecewise linear approx R=6}}    
    \label{fig:SoC_r6}
\end{subfigure}
\caption[]
{\small Time-SoC function and linear/piecewise linear approximation with various $R$ value} 
\label{fig:SoC_approx}
\end{figure*}

There exist multiple SoC charging functions of different batteries under different charging techniques. For example, a detailed time-SoC function provided by TeoTeslaFan in the Teals forum is publicly available \citep{TeslaForum}. This function indicates that for Tesla Model S, the battery could be fully charged in 110 minutes with the Tesla Supercharger. For our problem, the EV is assumed to be a smaller electric van whose charging time is less than 90 minutes, according to Alke's official website, instead of a normal private electric vehicle. Thus, some necessary changes are made to the function posted in \citep{TeslaForum} to ensure that the charging time from zero to full SoC is 90 minutes. This concave SoC function is estimated with linear approximation or piecewise linear approximation with a different number of segments. The representation of the actual time-SoC function and its approximation is shown in Figure \ref{fig:SoC_approx}, where the value $R$ is the number of segments used in piecewise linear approximation.

\subsection{Performance comparison between ALNS and VNS}
\begin{figure}[H]
    \centering
    \includegraphics[width = .6\textwidth]{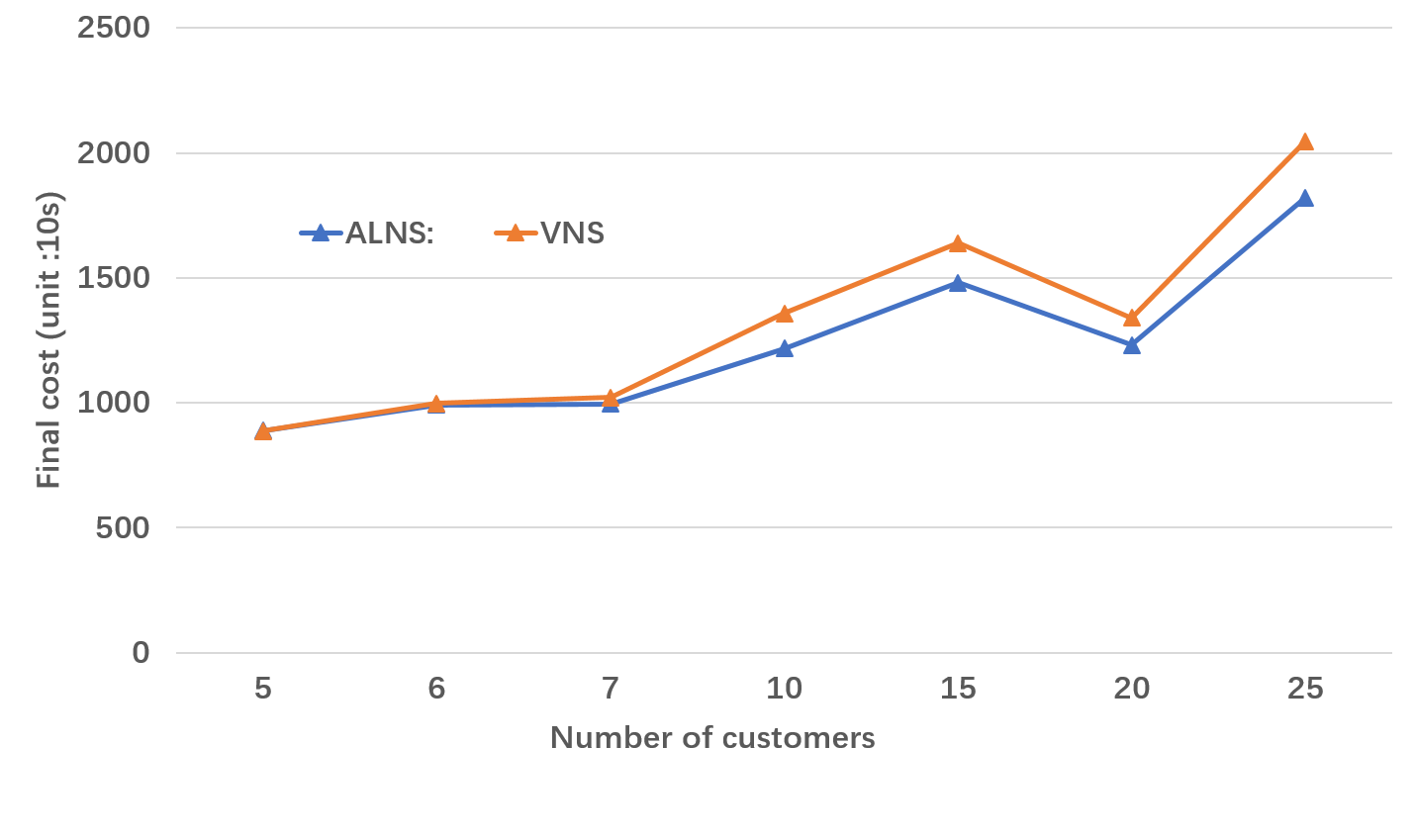}
    \caption{Comparison between ALNS and VNS }
    \label{fig:alns_vs_vns}
\end{figure}

First, the performance of the ALNS is compared to the VNS presented in \citep{zhu2022full}. The test results are shown in Figure \ref{fig:alns_vs_vns}, where the shown numbers are the average value of 20 randomly generated instances.

As can be seen, when the size of the tested instances is small ($|C| <= 10$), both heuristics have similar performance. However, as the instance size gets bigger ($|C| > 10$), ALNS continuously outperforms VNS, and the gap in the final cost is about 15\%. Note that we use the same computational time for both methods in this comparison. These test results indicate that ALNS performs better when dealing with EVTSPD-P instances. 
As explained earlier, the main advantage of ALNS over VNS is its ability to maintain the solution's feasibility. During the search process, the author found that many obtained solutions in VNS are discarded because the EV's driving range constraints are not satisfied. On the contrary, as ALNS has less number of iterations, the number of discarded solutions is much smaller. 

\subsection{A comparison between MILP formulation and ALNS method}
Here we analyze the performance of the ALNS heuristic by comparing it with solving various EVTSPD-P instances exactly using a commercial solver (CPLEX). When solving the EVTSPD-PN cases, we approximate the time-SoC function with the piecewise linear function of two segments ($|R|=2$). Test results on base EVTSPD-P instances is shown in Table \ref{tab:alns_performance1}.

Besides, we also conduct additional tests on EVTSPD-P with additional side constraints. The following variants are considered and the results are shown in Table \ref{tab:alns_performance2}.
\begin{itemize}
    \item \textit{Base}: the base case of EVTSPD-P without any side constraints
    \item \textit{LRT}: the EVTSPD-P variant considers the launch/retrieve time
    \item \textit{Range}: the EVTSPD-P variant which assumes the drone's flight range is dependent on the weight of parcel it carries.
    \item \textit{Loop}: the EVTSPD-P variant which allows the self-loop
\end{itemize}

\begin{landscape}
\begin{table}[H]
\begin{center}
\caption{Performance comparison between ALNS and CPLEX on base instances} \label{tab:alns_performance1}
\begin{tabular}{m{0.4cm} m{0.4cm} m{0.4cm} m{0.6cm} m{1cm} m{1cm} m{1cm} m{1cm} m{1cm} m{1cm} m{1cm} m{1cm} m{1cm} m{1cm} m{1cm}}
\toprule
\multicolumn{1}{l}{} & \multicolumn{1}{l}{} & \multicolumn{1}{l}{} & \multicolumn{1}{l}{} & \multicolumn{1}{l}{} & \multicolumn{4}{c}{Opt} & \multicolumn{4}{c}{RunningTime} & \multicolumn{2}{c}{Gap(\%)} \\
\cline{6-9}   \cline{10-13}    \cline{14-15}  \\
N & n & $n_{cs}$ & $\alpha$ & No.ins  & PL & PN & PL & PN  & PL & PN & PL & PN & PL & PN \\
\hline
8 & 4 & 2 & 1.5 & 20 & 796.58 & 772.20 & 796.58 & 792.20 & 17.20 & 10.80 & 5 & 5 & 0.00 & 2.59 \\
 & 4 & 2 & 2 & 20 & 743.25 & 741.24 & 743.25 & 741.24 & 17.20 & 10.80 & 5 & 5 & 0.00 & 0.00 \\
 & 4 & 2 & 2.5 & 20 & 705.69 & 699.21 & 705.69 & 699.21 & 17.10 & 10.80 & 5 & 5 & 0.00 & 0.00 \\
9 & 5 & 2 & 1.5 & 20 & 959.54 & 930.61 & 969.26 & 951.07 & 139.83 & 113.25 & 5 & 5 & 1.01 & 2.20 \\
 & 5 & 2 & 2 & 20 & 902.21 & 867.30 & 907.65 & 870.78 & 140.41 & 114.21 & 5 & 5 & 0.60 & 0.40 \\
 & 5 & 2 & 2.5 & 20 & 848.14 & 824.64 & 849.61 & 825.89 & 139.00 & 114.30 & 5 & 5 & 0.17 & 0.15 \\
10 & 6 & 2 & 1.5 & 10 & 1156.80 & 1122.78 & 1197.47 & 1157.64 & 2259.33 & 2804.00 & 20 & 20 & 3.52 & 3.10 \\
 & 6 & 2 & 2 & 10 & 1052.39 & 1019.85 & 1084.68 & 1040.02 & 2305.68 & 2814.27 & 20 & 20 & 3.07 & 1.98 \\
 & 6 & 2 & 2.5 & 10 & 980.32 & 960.29 & 1000.26 & 990.95 & 2345.17 & 2845.67 & 20 & 20 & 2.03 & 3.19 \\
14 & 8 & 3 & 1.5 & 10 & N/A & N/A & 1349.78 & 1295.78 & 3600 & 3600 & 20 & 20 & N/A & N/A \\
 & 8 & 3 & 2 & 10 & N/A & N/A & 1186.20 & 1136.56 & 3600 & 3600 & 20 & 20 & N/A & N/A \\
 & 8 & 3 & 2.5 & 10 & N/A & N/A & 921.53 & 902.60 & 3600 & 3600 & 20 & 20 & N/A & N/A \\
20 & 10 & 5 & 1.5 & 10 & N/A & N/A & 1657.97 & 1484.31 & 3600 & 3600 & 20 & 20 & N/A & N/A \\
 & 10 & 5 & 2 & 10 & N/A & N/A & 1431.35 & 1318.85 & 3600 & 3600 & 20 & 20 & N/A & N/A \\
 & 10 & 5 & 2.5 & 10 & N/A & N/A & 1237.06 & 1112.13 & 3600 & 3600 & 20 & 20 & N/A & N/A \\
36 & 20 & 8 & 1.5 & 10 & N/A & N/A & 1898.38 & 1634.4 & 3600 & 3600 & 20 & 20 & N/A & N/A \\
 & 20 & 8 & 2 & 10 & N/A & N/A & 1675.70 & 1424.23 & 3600 & 3600 & 20 & 20 & N/A & N/A \\
 & 20 & 8 & 2.5 & 10 & N/A & N/A & 1528.29 & 1396.48 & 3600 & 3600 & 20 & 20 & N/A & N/A \\
\bottomrule
\end{tabular}
\end{center}
\end{table}
\end{landscape}

\begin{table}[H]
\begin{center}
\caption{Performance comparison between ALNS and CPLEX on instances with side constraints} \label{tab:alns_performance2}
\begin{tabular}{m{0.4cm} m{0.4cm} m{0.4cm} m{0.4cm} m{0.6cm} m{0.9cm} m{0.9cm} m{0.9cm} m{0.9cm} m{0.9cm} m{0.9cm} m{0.9cm} m{0.9cm}}
\toprule
 & & & & &\multicolumn{2}{c}{\textbf{MaxLeg}} & \multicolumn{2}{c}{\textbf{LRT}} & \multicolumn{2}{c}{\textbf{Range}} &
  \multicolumn{2}{c}{\textbf{Loop}} \\
\cline{6-7}   \cline{8-9}    \cline{10-11}  \cline{12-13} \\
$|N|$ & n & $n_{cs}$ & $\alpha$ & No.ins & PL   & PN    & PL   & PN & PL   & PN & PL   & PN\\
\hline 
8 & 4 & 2 & 1.5 & 20 & 0.01 & 0.00 & 0.00 & 0.00 & 0.00 & 0.00 & N/A & N/A \\
 & 4 & 2 & 2 & 20 & 0.00 & 0.00 & 0.00 & 0.00 & 0.00 & 0.00 & N/A & N/A \\
 & 4 & 2 & 2.5 & 20 & 0.00 & 0.00 & 0.00 & 0.00 & 0.00 & 0.00 & N/A & N/A \\
9 & 5 & 2 & 1.5 & 20 & 1.02 & 1.45 & 1.75 & 1.62 & 1.32 & 1.04 & N/A & N/A \\
 & 5 & 2 & 2 & 20 & 1.69 & 1.24 & 1.56 & 1.46 & 0.79 & 0.78 & N/A & N/A \\
 & 5 & 2 & 2.5 & 20 & 1.30 & 1.41 & 1.32 & 1.48 & 1.74 & 1.35 & N/A & N/A \\
10 & 6 & 2 & 1.5 & 10 & 3.05 & 3.07 & 2.98 & 3.45 & 3.14 & 3.21 & N/A & N/A \\
 & 6 & 2 & 2 & 10 & 3.04 & 3.14 & 3.62 & 3.11 & 3.12 & 3.23 & N/A & N/A \\
 & 6 & 2 & 2.5 & 10 & 2.97 & 3.25 & 3.75 & 3.18 & 3.29 & 3.37 & N/A & N/A \\
\bottomrule
\end{tabular}
\end{center}
\end{table}

As can be seen from Table \ref{tab:alns_performance1}, for the base cases, the commercial solver can only solve instances containing six customers and two charging stations within one hour. Note that for instances of this size, the actual number of nodes in the augmented network is 10. For these instances, the average optimality gap of ALNS is less than 3.5 \% when the computational time of ALNS is set to five seconds. For instances containing more than ten nodes in the augmented network, CPLEX fails to converge in one hour.\par

A comparison between the final cost of EVTSPD-PL and EVTSPD-PN instances indicates that, on average, the latter has about 7\% lower cost than the former case, and this gap seems to grow as the size of the instances increases. This result indicates that using a piecewise linear approximation of the time-SoC function can render higher quality results than a linear approximation. The drone's speed is a significant factor that affects the final cost, as the final delivery time of EVTSPD-P with $\alpha = 2.5$ is about 12\% lower than that with $\alpha = 1.5$.\par

The test results shown in Table \ref{tab:alns_performance2} indicate that the optimality gap for instances with additional side constraints is similar to that of the base instances for MaxLeg, LRT, and Range variants. However, for the \textit{Loop} variant, CPLEX fails to terminate in one hour, even for a small network with four customers and two charging stations. This is because additional dummy copies of the customer nodes are needed for the MILP formulation to handle the \textit{Loop} variant.

\subsection{A comparison of EVTSPD-P with different $R$ value}

In this subsection, we compare the final cost of EVTSPD-P when using different time-SoC function approximations. The actual time-SoC function and its approximation with different $|R|$ values are shown in Figure \ref{fig:SoC_approx}. Obviously, as the $|R|$ value increases, the accuracy of the approximation also increases. We adopt a conservative method for all the linear/piecewise linear approximations to ensure the obtained solution is feasible when using the actual time-SoC function.\par 
We first conduct the tests on randomly generated instances, and the final delivery time with different $|R|$ values are shown in Table \ref{tab:various_R}. As expected, the average final cost of EVTSPD-P instances decreases as the $|R|$ value increases. The average final cost with $|R|=6$ is about 15728.3 seconds, while that with linear approximation is about 17643.2 seconds, 10.88\% higher than when $|R|=6$. These results indicate that using linear time-SoC function approximation might generate a sub-optimal solution compared to using piecewise linear function approximation with a high $|R|$ value.\par

\begin{figure}[H]
    \centering
    \includegraphics[width = .6\textwidth]{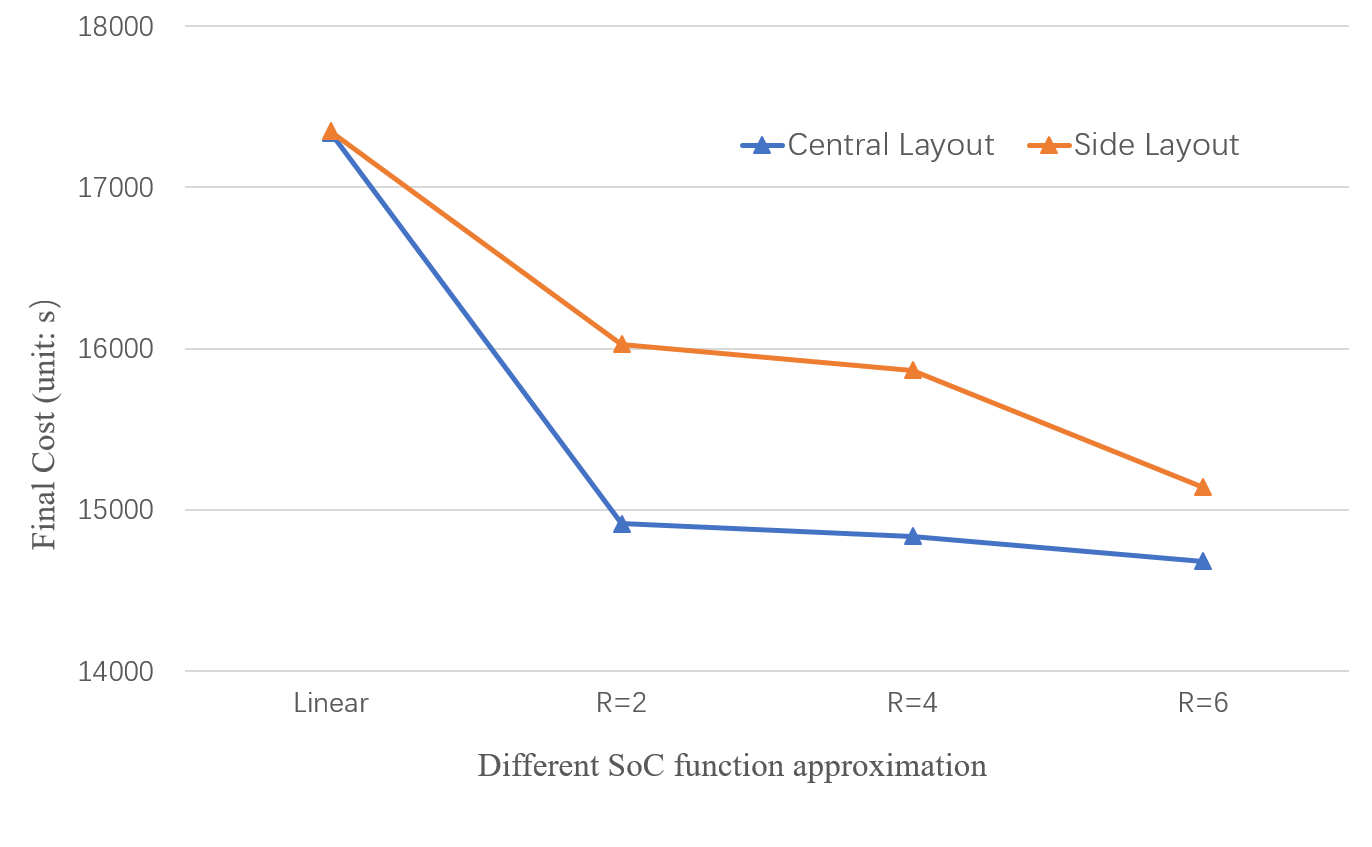}
    \caption{Final cost of Austin network with various $|R|$ value}
    \label{fig:AustinNet_variousR}
\end{figure}

Besides, we also test the final delivery time of the downtown Austin network described in \citep{zhu2022full} with various $|R|$ values. The results are shown in Figure \ref{fig:AustinNet_variousR}, which has similar results as described above. For both network layouts, the final delivery time decreases gradually as $|R|$ increases, and the gap between $|R|=1$ and $|R|=6$ is about 15.0\% and 12.7\% for central and side layouts, respectively.  

\begin{table}[H]
\begin{center}
\caption{Performance comparison of EVTSPD with different $|R|$ value} \label{tab:various_R}
\begin{tabular}{llllllll}
\toprule
 &  &  &  & \multicolumn{4}{c}{Final cost(unit: 10s)} \\
 \cline{5-8}
N & C & N\_CS & No.Ins & PL & R=2 & R=4 & R=6 \\
\hline
10 & 6 & 2 & 10 & 1131.52 & 1089.01 & 1084.82 & 1084.09 \\
20 & 10 & 5 & 10 & 1630.02 & 1537.38 & 1507.80 & 1494.32 \\
30 & 16 & 7 & 10 & 1815.57 & 1712.60 & 1677.74 & 1664.02 \\
40 & 20 & 10 & 10 & 1997.34 & 1757.35 & 1697.93 & 1690.60 \\
50 & 26 & 12 & 10 & 2247.16 & 1935.02 & 1934.24 & 1931.11 \\
\textbf{Average} &  &  &  & \textbf{1764.32} & \textbf{1606.27} & \textbf{1580.51} & \textbf{1572.83} \\
\bottomrule
\end{tabular}
\end{center}
\end{table}

\section{Conclusions and future work}

As an extension of the work in \citep{zhu2022full}, in this paper, we address another EVTSPD variant which assumes that the EV is a plug-in hybrid electric vehicle that could be partially charged in the charging stations which are assumed to be connected to a plug-in electric grid. With this assumption, the electric vehicle in EVTSPD can be charged at CS nodes with any amount of energy, while the required charging time depends on the initial state-of-charge level before charging. In this way, both the amount of charge energy and the charging time at each station node are decision variables. This problem, which aligns with the real-world application, is called EVTSPD-P, where the letter "P" stands for partial recharging. \par

It turns out that adopting the partial-recharge assumption indicates that the pre-constructed multigraph introduced in \citep{zhu2022full} is not applicable in EVTSPD-P, as the charging time at each CS node is not fixed. Thus, the charging station nodes have to be included in the network. Furthermore, "dummy" copies of the CS nodes also need to be included to enable potential multiple visits to the same charging station. For these reasons, the resulting network, which is termed as "the augmented network" in this paper, is much bigger in size than the multigraph. A MILP formulation defined in this augmented network is proposed in this paper to solve EVTSPD-P. \par

Apparently, one key step in modelling EVTSPD-P is to account for the time-SoC function, which is a concave function that specifies the relationship between the charging time and state-of-charge level. To incorporate this function into the model, two different approximation techniques, namely, linear function approximation and piecewise linear function approximation, are adopted, resulting in two different variants of EVTSPD-P, namely, EVTSPD-PL and EVTSPD-PP, where "PL" stands for "partial linear" while "PP" stands for "partial piecewise linear". Obviously, EVTSPD-PP is more accurate than EVTSPD-PL, with the additional complexity of modeling extra decision variables. Modelling EVTSPD-PL is straightforward with the proposed MILP formulation, while the latest technique proposed in \citep{zuo2019new} is used to model EVTSPD-PP.\par

Considering the enormous number of constraints involved in EVTSPD-PL and EVTSPD-PP, solving the MILP formulation via the CPLEX solver only works for small-size instances. To handle EVTSPD-P instances with practical size, we propose a specially designed adaptive neighborhood search meta-heuristic method. LNS is similar to VNS that is described in \citep{zhu2022full}, in that both approaches are neighborhood-search-based methods. However, VNS only adopts several relatively simply neighborhood structures while LNS utilizes a series of "heavy" methods when searching for a new feasible solution. VNS is suitable for solving problems and it is easy to identify new feasible solutions. However, for EVTSPD-P, with the extra EV's driving range constraints, it is difficult for VNS to maintain all the new solution's feasibility with simply neighborhood structures. Instead, ALNS has shown its reasonable performance in handling problems with tight constraints. Within each iteration of LNS, the quality of new solutions found is better than VNS as the former approach adopts a more complex neighborhood structure for finding new solution, and it is relatively easy for LNS to maintain the solution's feasibility. In our paper, we adopt an adaptive LNS, which indicates that the performance of each destroy and repair method is evaluated and updated throughout the search process. Besides, we propose a new repair method specifically designed for EVTSPD-P, which is based on constraint programming. 
 \par

In the numerical analysis section, a thorough comparison between the MILP formulation, ALNS, and VNS described in \citep{zhu2022full} is conducted. The results indicate that ALNS outperforms VNS when the number of nodes in the augmented network exceeds ten and that the ALNS has an average optimality gap of less than 3.5\%. Besides, test results on different approximation parameters ($|R|$ value) show that the piecewise function with six segments has an average of 10\% less cost than linear function approximation. \par

With the current research, several problems remain unsolved and might be a good research direction for the future. First of all, it is still a challenge to obtain the exact solution of EVTSPD-P instances with practical size ($|C|>20$). Secondly, the future search could consider uncertainties and randomness of the problem, such as the customers' demands or the flight time/speed of the drone.

\section{Acknowledgements}
This research is based on work supported by the National Science Foundation under Grant No. 1826230, 1562109/1562291, 1562109, 1826337, 1636154, and 1254921. This work is also supported by the Center for Advanced Multimodal Mobility Solutions and Education (CAMMSE).

\appendix
\section{MIP formulation of CP repair method}
This appendix describes the mixed-integer program sub-problem for the CP repair method to determine the drone sorties. This method takes the partial solution and a set of selected nodes as input and seeks to serve these nodes by drone, inserting these sorties into the truck's route. At this point in the process, we have eliminated all  existing drone sorties, so the partial solution only contains several trucks' routes. This formulation is initially proposed in \citep{yurek2018decomposition}. The definition of the indexes, sets, decision variables, and parameters are described below:

\subsection{Indexes}
{\renewcommand\arraystretch{1.0}
\noindent\begin{longtable*}
{@{}l @{\quad:\quad} p{15cm}@{}}
% {@{}>{\raggedright}p{2cm}lp{7.3cm}@{}}
$i,j$ & node\\
$k$ & position \\
$p$ & drone sortie \\
\end{longtable*}}

\subsection{Sets}
{\renewcommand\arraystretch{1.0}
\noindent\begin{longtable*}
{@{}l @{\quad:\quad} p{15cm}@{}}
% {@{}>{\raggedright}p{2cm}lp{7.3cm}@{}}
$C$ & set of all customer nodes\\
$D$ & set of customers that would be served by the drone \\
$S$ & set of all potential drone sorties \\
\end{longtable*}}

\subsection{Parameters}
{\renewcommand\arraystretch{1.0}
\noindent\begin{longtable*}
{@{}l @{\quad:\quad} p{15cm}@{}}
% {@{}>{\raggedright}p{2cm}lp{7.3cm}@{}}
$d_{s}$ & duration of sortie $s$\\
$f_{is}$ & binary parameter which equals 1 if sortie $s$ starts from node $i$, and 0 otherwise \\
$a_{is}$ & binary parameter which equals 1 if sortie $s$ serves node $i$, and 0 otherwise \\
$l_{is}$ & binary parameter which equals 1 if sortie $s$ ends at node $i$, and 0 otherwise \\
$N$ & number of customer nodes in the instance \\
$t_{i}$ & arrival time of the truck at node $i$  \\
$m_{k}$ & node assigned to position $k$ in the truck's route  \\
$r_{t}$ & current truck's route \\
$\{0, C+1\}$ & the depot set \\
\end{longtable*}}

\subsection{Decision variables}
{\renewcommand\arraystretch{1.0}
\noindent\begin{longtable*}
{@{}l @{\quad:\quad} p{15cm}@{}}
% {@{}>{\raggedright}p{2cm}lp{7.3cm}@{}}
$x_{s}$ & binary parameter which equals 1 if sortie $s$ is chosen in the final solution, and 0 otherwise\\
$w_{i}$ & waiting time of the truck at node $i$ \\
\end{longtable*}}

\vspace{0.2cm} 
\setcounter{equation}{0}
\noindent Objective :
\begin{align}
\min \quad
& \sum_{j \in C\setminus D \cup \{C+1\} }w_{i}
\end{align}

\vspace{0.5cm}
\noindent The objective is to minimize the total waiting time that is incurred by adding drone sorties to the current truck's route. \par

\hfill \break
\noindent \textbf{s.t.}
\begin{align}
\sum_{s \in S}\{d_{s}l_{is}x_{s}-(t_{i}l_{is}x_{s}-\sum_{j \in C\setminus D \cup \{0\}} t_{j}f_{js}x_{s}) \} & \leq w_{i} && \forall i \in C\setminus D \cup \{C+1\} \\
\sum_{s\in S} a_{is}x_{s}& = 1  && \forall i \in D\\
\sum_{s\in S} f_{is}x_{s}& \leq 1  && \forall i \in C\setminus D \cup \{0\}\\
\sum_{s\in S} l_{is}x_{s}& \leq 1  && \forall i \in C\setminus D \cup \{C+1\}\\
\sum_{s\in S} f_{m_{k}s}x_{s} + \sum_{s\in S} l_{m_{k}s}x_{s} - 2\left( 1- \sum_{s\in S} f_{m_{i}s}l_{m_{j}s}x_{s}\right) & \leq 0  && \forall i=0,1,..n-1, j = i+2,...,C+1, i \leq k \leq j, i\neq j\\
x_{s} &\in \{0,1\} && s \in S \\
w_{i} &\geq 0 && i\in C\setminus D \cup \{C+1\}
\end{align}

Constraint (2) specifies that the waiting time is incurred when the truck waits for the drone. Constraints (3), (4), and (5) indicate that every node can only be a launch node, drone node, and retrieve node for at most one sortie, respectively. Constraint (6) indicates that the drone cannot be re-launched before it is retrieved. Constraints (7) and (8) specify the domain of the decision variables.

\bibliography{library}
\end{document}